\documentclass[12pt,psamsfonts]{amsart}
\usepackage[left=3cm,top=2.5cm,right=3cm,bottom=2.5cm]{geometry}

\newtheorem {teo} {Theorem} [section]
\newtheorem {propo} [teo] {Proposition}
\newtheorem {corol} [teo] {Corollary}

\usepackage{amsmath,amssymb}
\usepackage{graphicx}
\usepackage[dvips]{psfrag}
\usepackage{color}
\usepackage{psfrag} % See the section on figures.
\usepackage{overpic} % See the section on figures.
\usepackage{doi}
\usepackage{subcaption}

\graphicspath{{figs/}}
 
\newtheorem{theorem}{Theorem}
\newtheorem{proposition}[theorem]{Proposition}

\newtheorem{lemma}[theorem]{Lemma}
\newtheorem {remark}[theorem]{Remark}
\newtheorem{definition}[theorem]{Definition}

\newcommand{\w}{{\omega}}
\newcommand{\te}{{\theta}}
\newcommand{\wz}{{\omega_0}}
\newcommand{\rz}{{\rho_0}}

\title[]{Cyclicity of  Centers on Center Manifolds in a 3D Chaotic System with a Four-Wing Attractor}
\author[V. Gusson and C. Pessoa]{}

  \subjclass[2020]{34C07, 34C25}
   \keywords{Cyclicity; Center-Focus problem; Isochronicity; 3D chaotic systems}

\begin{document}
 \maketitle

\centerline{\scshape Vitor Gusson and Claudio Pessoa}
\medskip

{\footnotesize \centerline{Universidade Estadual Paulista (UNESP),} \centerline{Instituto de Bioci\^encias Letras e Ci\^encias Exatas,} \centerline{R. Cristov\~ao Colombo, 2265, 15.054-000, S. J. Rio Preto, SP, Brazil }
\centerline{\email{vitor.gusson@unesp.br},\;\email{c.pessoa@unesp.br}}

\medskip

\bigskip

\begin{quote}{\normalfont\fontsize{8}{10}\selectfont
{\bfseries Abstract.} In this work, we investigate the conditions that guarantee the existence of centers on the center manifold, arising from Hopf points, in the new three-dimensional quadratic chaotic system introduced by B. Khaled et al. in 2024 in the {\it Int. J. Data Netw. Sci}. For some of the system’s Hopf points, we solve the center-focus problem on the center manifold, analyzing both its isochronicity and cyclicity. Our results significantly improve the previously known lower bound on the number of limit cycles bifurcating from Hopf points in this system, as established by B. M. Mohammed in 2025 in the {\it Int. J. Bifurc. Chaos Appl. Sci. Eng.}.
\par}
\end{quote}

\section{Introduction and Main Results}

One of the main interests of the field of dynamical systems and the qualitative theory of ordinary differential equations is the study of chaotic systems. Although behaviors sensitive to initial conditions were observed by Henri Poincaré in mathematical models of celestial mechanics, particularly in the formulation and study of the three-body problem (see \cite{Poincare}), it is Lorenz who is credited with the discovery of the first chaotic system. In his 1963 work, ``Deterministic Non-periodic Flows" (see \cite{Lorenz}), Edward Lorenz formulated a nonlinear system of three equations to model atmospheric phenomena in meteorology and systematically presented the incidence of chaos for certain parameter values. The Lorenz system is given by
\begin{equation*}
	\label{Lorenz}
	\begin{array}{l}
		\dot{u} =\sigma(v-u), \\
		\dot{v} = \rho u-v-u w, \\
		\dot{w} = u v-\beta w, \\
	\end{array}
\end{equation*}
where $\sigma,\rho,\beta>0$.

Since then, the field of chaotic systems has expanded considerably. In an attempt to better understand and characterize the problem, research has focused on identifying chaos using tools such as Lyapunov exponents (\cite{LyapExp68,LyapExp16,LyapExp85}) and bifurcation diagrams (\cite{BifDiag78,BifDiag04}), which serve, respectively, to measure the exponential divergence rate between close trajectories and to analyze the change in qualitative behavior as parameters vary. Furthermore, there is interest in understanding the geometry of strange attractors---the sets to which the system’s trajectories converge---(see \cite{StrgAttr76,Lorenz,StrgAttr71}), as well as in the control (see \cite{Ctrl06,Ctrl90,Ctrl92}) and synchronization of chaos (see \cite{Sync02,Sync90}), since these phenomena appear in a wide range of physical, biological, and economic models (see \cite{Strogatz15}, for example).

Another important aspect of the qualitative theory of ordinary differential equations is the study of periodic orbits. This includes problems such as determining center conditions for centers on center manifolds of three-dimensional differential systems, analyzing the occurrence of limit cycles due to Hopf bifurcations, and establishing upper or lower bounds on the number of limit cycles. The latter is related to Hilbert’s 16th problem (see \cite{16thHilbert}). Several well-known three-dimensional chaotic systems have already been studied from these perspectives. In \cite{Chen03,LlibrePessoa15,MelloCoelho09}, the authors investigated the occurrence of Hopf bifurcations in the Chen (see \cite{Chen99}), Shimizu–Morioka (see \cite{ShimizuMorioka80}), and Lü (see \cite{Lu02}) systems, respectively. In \cite{Wang12}, multiple limit cycles arising from Hopf bifurcations in the Chen system were analyzed. In \cite{MahdiPessoa11}, conditions for the existence of centers on the center manifold of the Lü system were established.
In \cite{GouveiaQueiroz24}, the authors determined the center conditions for the Lorenz and Rössler systems (see \cite{Rossler77}) and provided a lower bound on the number of limit cycles that can arise when perturbing these centers. Other works on lower bounds for the number of limit cycles that bifurcate from perturbed centers on center manifolds of three-dimensional differential systems can be found in \cite{Gar18,Ke25,Pessoa24,SanchezTorre22,Tian14, Yu15}.

In the pursuit of a better understanding of chaotic models and their various characteristics, there has been a growing interest in the study of chaotic systems with multi-wing attractors. These systems show promising potential in areas such as secure communications and image encryption. In references \cite{OnetoFourWing16,TwoWing20, OnetoFourWing18}, for example, the authors find systems with one to four wing attractors and a systematic study of the occurrence of chaos. In \cite{Mohammed25}, a chaotic system with four-wing attractors was considered, which was originally introduced and studied in \cite{FourWing24}. In that work, the author analyzed the occurrence of Hopf bifurcations around equilibrium points and investigated the system's behavior at infinity. 

In this paper, we aim to deepen the understanding of the local dynamics around the equilibrium points of the chaotic system with four-wing attractors and to improve upon some of the results presented in \cite{Mohammed25}. Specifically, we address the center-focus problem on center manifolds in cases where the system’s singular points are Hopf points. Additionally, we investigate the isochronicity of the classified centers and establish lower bounds on the number of limit cycles that bifurcate from these centers.

As in \cite{Mohammed25}, we will consider the system
\begin{equation}\label{sist3Dchaotic}
	\begin{array}{l}
		\dot{x}=a(y-x)+y z,\\
		\dot{y}=b x +c y-x z,\\
		\dot{z}=-d z+x y+1.
	\end{array}
\end{equation}
Let $\Delta=(a+b)^2+4ac$. When $d\neq0$, this system has the following singular points:
\begin{align*}
	E_1&=\left(0,0,\dfrac{1}{d}\right); \\
	E_2^{\pm}&=\Bigg( \pm\dfrac{\sqrt{-1+b d+2 c d-\frac{-b+\sqrt{\Delta}\;(-1+b d)}{a}}}{\sqrt{2}},\\
	&\quad \mp\dfrac{\left(a+b+\sqrt{\Delta}\right)\sqrt{-1+b d+2 c d-\frac{-b+\sqrt{\Delta}\;(-1+b d)}{a}}}{2\sqrt{2}c},\dfrac{1}{2} \left(-a + b - \sqrt{\Delta}\right)\Bigg); \\
	E_3^{\pm}&=\Bigg( \pm\dfrac{\sqrt{-1+b d+2 c d+\frac{b+\sqrt{\Delta}\;(-1+b d)}{a}}}{\sqrt{2}},\\
	&\quad \mp\dfrac{\left(a+b-\sqrt{\Delta}\right)\sqrt{-1+b d+2 c d+\frac{b+\sqrt{\Delta}\;(-1+b d)}{a}}}{2\sqrt{2}c},\dfrac{1}{2} \left(-a + b + \sqrt{\Delta}\right)\Bigg).
\end{align*}
Note that, point $E_1$ always exists. Assuming $\Delta\geq0$ and $a c\neq0$, then $E_2^\pm$ and $E_3^\pm$ exist for $-1+b d+2 c d-\frac{-b+\sqrt{\Delta}\;(-1+b d)}{a}\geq0$ and $-1+b d+2 c d+\frac{b+\sqrt{\Delta}\;(-1+b d)}{a}\geq0$, respectively.
Now, for $d=0$, we will consider the following parameters sets:
\begin{align*}
	W_1&=\{(a,b,c)\in\mathbb{R}^3\;:\;a>0,\;c\neq0,\;b<-a-\sqrt{\Delta}\}, \\
	W_2&=\{(a,b,c)\in\mathbb{R}^3\;:\;a<0,\;c\neq0,\;b>-a-\sqrt{\Delta}\}, \\	W_3&=\{(a,b,c)\in\mathbb{R}^3\;:\;a>0,\;c\neq0,\;b<-a+\sqrt{\Delta}\}, \\	W_4&=\{(a,b,c)\in\mathbb{R}^3\;:\;a<0,\;c\neq0,\;b>-a+\sqrt{\Delta}\},
\end{align*}
with $\Delta>0$, and the points

\begin{equation}\label{ptsE4E5}
	\begin{array}{r@{\,=\,}l}
		E_4^\pm & \Bigg( \pm\dfrac{\sqrt{-\dfrac{a+b+\sqrt{\Delta}}{a}}}{\sqrt{2}},\mp\dfrac{\left(a+b-\sqrt{\Delta}\right)\sqrt{-\frac{a+b+\sqrt{\Delta}}{a}}}{2\sqrt{2}c}, \dfrac{1}{2} \left(-a + b + \sqrt{\Delta}\right)\Bigg);\\
		E_5^\pm & \Bigg( \pm\dfrac{\sqrt{-\dfrac{a+b-\sqrt{\Delta}}{a}}}{\sqrt{2}},\mp\dfrac{\left(a+b+\sqrt{\Delta}\right)\sqrt{-\frac{a+b-\sqrt{\Delta}}{a}}}{2\sqrt{2}c}, \dfrac{1}{2} \left(-a + b - \sqrt{\Delta}\right)\Bigg).
	\end{array}
\end{equation}

Hence, if $(a,b,c)\in W_1\cup W_4$, points $E_4^\pm$ and $E_5^\pm$ are the only singular points of system \eqref{sist3Dchaotic}. If $(a,b,c)\in W_2$, we just have singular points $E_4^\pm$ and, if $(a,b,c)\in W_3$, only singular points $E_5^\pm$. Note, also, that system \eqref{sist3Dchaotic} is invariant under the change of coordinates $(x,y,z)\mapsto(-x,-y,z)$.

We have the following new results on the investigation of the center-focus problem on center manifolds for system \eqref{sist3Dchaotic}.

\begin{theorem}\label{Teo1}
	Consider system \eqref{sist3Dchaotic} with $d\ne0$. If \( a = c \) and $(1+cd)(1-bd)- c^2 d^2 > 0$, then singular point $E_1$ is a Hopf point of the system. Furthermore, under these conditions, the Hopf point  $E_1$ is a center on the center manifold of the system if and only if $b = c = 0$, and it is not an isochronous center.
\end{theorem}

\begin{theorem}\label{Teo2}
	Consider system \eqref{sist3Dchaotic} with $d=0$. If $a=-c$, $c>0$ and $b>3c$ (resp. $a=-c$, $c<0$ and $b<3c$), then singular points $E_4^\pm$ (resp. $E_5^\pm$) are a pair of Hopf points of system \eqref{sist3Dchaotic}, which are foci on its center manifold.
\end{theorem}

%\begin{theorem}\label{Teo3}
%	Consider system \eqref{sist3Dchaotic} with $d=0$. Then, points $E_5^\pm$ are a pair of Hopf points for conditions $a=-c$, $c<0$ and $b<3c$. Also, under these conditions, $E_5^\pm$ are a focus on the center manifold.
%\end{theorem}

Once the centers of system \eqref{sist3Dchaotic} have been classified, a natural question arises: how many limit cycles can bifurcate from its centers under small perturbations of the parameters? We provide a partial answer to this question: we have shown that at least three limit cycles can bifurcate from the center described in Theorem \ref{Teo1}. This result improves upon that of \cite{Mohammed25}, where it was shown that at least one limit cycle bifurcates from the singular point $E_1$ via a Hopf bifurcation. Furthermore, when considering quadratic perturbations of this center, we have shown that at least five limit cycles can bifurcate. Our results are summarized in follows theorems. 

\begin{theorem}\label{Teo4}
	There exist values of the parameters $(a, b, c, d)$ in system \eqref{sist3Dchaotic} for which the system has at least three limit cycles, bifurcating from the singular point $E_1$ under the center conditions described in Theorem \ref{Teo1}.
\end{theorem}

\begin{theorem}\label{Teo5}
	Consider system \eqref{sist3Dchaotic}, with the center conditions established in Theorem \ref{Teo1}, then at least five limit cycles can bifurcate from center $E_1$ through quadratic perturbations of the system.
\end{theorem}

The paper is organized as follows: In Section \ref{secCenter}, we give more details on the center-focus problem on center manifolds, i.e., we present an algorithm for determining whether a Hopf point is a center or a focus for the system restricted to its center manifold. In Section \ref{secIsoc}, we briefly introduce the issue of isochronicity for centers and provide an algorithm that determines when a center on the center manifold is isochronous. In Section \ref{secProof}, we prove Theorems \ref{Teo1} and \ref{Teo2}. In Section \ref{secSimul}, we do some numerical simulations that illustrate the dynamic behavior near points $E_1$, $E_4^\pm$ and $E_5^\pm$, under the conditions of Theorems \ref{Teo1} and \ref{Teo2}. Finally, in Section \ref{secCyc}, we show some results to study the cyclicity of three-dimensional differential systems and perform the proofs of Theorems \ref{Teo4} and \ref{Teo5}.

\section{Center-Focus Problem}\label{secCenter}
In this section, we present the center-focus problem for three-dimensional differential systems that have a monodromic singular point of the Hopf type, that is, singular points for which the associated eigenvalues are a pair of purely imaginary eigenvalues and a non-zero real eigenvalue. For these types of singular points, we present an algorithm, already known in the literature and found in detail in \cite{Edneral12} (see also \cite{Bibikov}), to distinguish a center from a focus.

A system $\dot{\mathbf{x}}=X(\mathbf{x})$, where $\mathbf{x} \in \mathbb{R}^3$, $X$ is an analytical vector field defined in an open set $U \subset \mathbb{R}^3$, and $X(0)=0$ such that $0$ is a Hopf point, can always be written in the following form, after a linear change of coordinates and a rescaling of time, if necessary:
\begin{equation}
	\label{formnorm}
	\begin{array}{l}
		\dot{u} = -v + P(u,v,w), \\
		\dot{v} = u + Q(u,v,w), \\
		\dot{w} = \lambda w + R(u,v,w), \\
	\end{array}
\end{equation}
with $\lambda \neq 0$ and $P, Q, R$ analytical functions with Taylor series expansions having only terms of order greater than or equal to $2$. System \eqref{formnorm} satisfies the Center Manifold Theorem, that is, for all $r \in \mathbb{N}$, there exists an invariant manifold $W^c$ of class $\mathcal{C}^r$ tangent to the $xy$-plane at the origin. The manifold $W^c$ is called the \textit{center manifold}, and a proof of this result is seen in \cite{Kel67}. Furthermore, since the origin of \eqref{formnorm} is a Hopf point, it can be either a center or a focus for the system restricted to a center manifold $W^c$. In three-dimensional differential systems, the center-focus problem fundamentally lies in deciding whether the singular point at the origin of system \eqref{formnorm}, restricted to a center manifold, is a center or a focus.

Before presenting the fundamental result on which the algorithm is based, let us recall some important basic definitions.

The vector field $\mathcal{X}$, associated with system \eqref{formnorm} and defined in the neighborhood of the origin, is represented by:
\begin{equation}\label{eqCampo}
	\mathcal{X}=(-v+P)\frac{\partial}{\partial u}+(u+Q)\frac{\partial}{\partial v}+(\lambda w+R)\frac{\partial}{\partial w}.
\end{equation}
A local analytic first integral of system \eqref{formnorm} corresponds to a non-constant differentiable function $\Psi$, mapping from a neighborhood of the origin in $\mathbb{R}^3$ to $\mathbb{R}$ and remaining constant along the system's trajectories, that is,
\begin{equation*}
	\mathcal{X} \Psi=(-v+P)\frac{\partial \Psi}{\partial u}+(u+Q)\frac{\partial\Psi}{\partial v}+(\lambda w+R)\frac{\partial\Psi}{\partial w}=0.
\end{equation*}
More generally, a formal first integral of system \eqref{formnorm} is defined as a non-constant formal power series $\Psi$ in variables $u$, $v$, and $w$, centered at the origin, which also satisfies the condition $\mathcal{X}\Psi=0$.

The following result, with its proof available in \cite[p. 2617]{Edneral12}, provides a theoretical solution in the study of the center-focus problem in three-dimensional differential systems. 

\begin{theorem}\label{TeoEdneral}
	Consider system \eqref{formnorm}. The following statements are equivalent.
	\begin{itemize}
		\item[$I$.] The origin is a center for the vector field restricted to a center manifold $\mathcal{X}\vert_{W^c}$;
		\item [$II$.] System \eqref{formnorm} admits a local analytic first integral that can be given as $\Psi(u,v,w)=u^2+v^2+\cdots$, where the dots are higher-order terms;
		\item [$III$.] System \eqref{formnorm} admits a formal first integral.
	\end{itemize}
\end{theorem}

Essentially, the algorithm to be detailed in this work is based on the equivalence $I \iff III$. This property will enable us to determine the necessary conditions for a Hopf point at the origin of system \eqref{formnorm} to constitute a center on the center manifold. However, we would also like to mention the work \cite{Wang10}, which employs formal series to address the center-focus problem and provides an alternative to the algorithm presented here. The latter has already been applied to chaotic systems, as in \cite{Wang14}.

We start this investigation by complexifying system \eqref{formnorm} through the complex variables $x=u+i v$, $y=\bar{x}=u-i v$, and $z=w$. We obtain the following differential system in its complex form:
\begin{align*}
	\dot{x} &= ix + X_1(x, y, z), \\
	\dot{y} &= -iy + X_2(x, y, z), \\
	\dot{z} &= \lambda z + X_3(x, y, z).
\end{align*}
Here, $X_i$ represents sums of homogeneous polynomials, with degrees ranging from $2$ to the degree $N$ of system \eqref{formnorm}. In fact, the preceding system can be expressed as follows:
\begin{equation} \label{formComplex}
	\begin{aligned}
		\dot{x} &= ix + \sum_{j+k+l=2}^{N} a_{jkl}x^j y^k z^l, \\
		\dot{y} &= -iy + \sum_{j+k+l=2}^{N} b_{jkl}x^j y^k z^l, \\
		\dot{z} &= \lambda z + \sum_{j+k+l=2}^{N} c_{jkl}x^j y^k z^l.
	\end{aligned}
\end{equation}
The coefficients $b_{jkl}$ are the conjugates of $a_{jkl}$, and $c_{jkl}$ are determined such that the sum $\sum_{j+k+l=2}^{N} c_{jkl}x^j y^k z^l$ yields a real value for all $x \in \mathbb{C}$ and all $z \in \mathbb{R}$. Our next step is to investigate the existence of a first integral $\Psi = \Psi(x, y, z)$ for system \eqref{formComplex}, which is equivalent to finding a first integral $\Phi = \Phi(u, v, w)$ for system \eqref{formnorm}. Note that the function $\Psi$, if it exists, can be expressed, by Theorem \ref{TeoEdneral}, in the following form:
\begin{equation}\label{eq:psi_form}
	\Psi(x, y, z) = xy + \sum_{j+k+l\ge3} d_{klj}x^k y^l z^j.
\end{equation}
Let us denote $\tilde{\mathcal{X}}$ as the vector field associated with system \eqref{formComplex} in $\mathbb{C}^3$, so it must satisfy the condition $\tilde{\mathcal{X}}\Psi=0$. 
Denoting by $L_{k_1 k_2 k_3}$ the coefficient of $x^{k_1}y^{k_2}z^{k_3}$ in $X\Psi$, except when $(k_1, k_2, k_3) = (k, k, 0)$ for a positive integer $k$, we can uniquely solve the equation $L_{k_1 k_2 k_3} = 0$ for $d_{k_1 k_2 k_3}$ in terms of the known quantities $d_{\alpha \beta \gamma}$ such that $\alpha + \beta + \gamma < k_1 + k_2 + k_3$. Hence, if $L_{kk0} = 0$ for all $k \in \mathbb{N}$, a formal first integral $\Psi$ exists. When the coefficient $L_{kk0}$ is nonzero, an obstruction to the existence of the formal series $\Psi$ arises. Then, we obtain 
\begin{equation}\label{eqFocal}
	\tilde{\mathcal{X}}\Psi(x, y, z) = L_{1}(xy)^2 + L_{2}(xy)^3 + \dots=\sum_{j=2}^{\infty} L_{j-1} (x y)^j,
\end{equation}
where $L_{k-1}=L_{kk0}$, for $k\geq2$.
The quantities $L_K$ are referred to as the {\it focus quantities} of system \eqref{formComplex} or {\it focus quantities}. The vanishing of all these focus quantities is a necessary and sufficient condition for the origin of system \eqref{formnorm} to be a center on the center manifold.

Another approach to studying the center-focus problem is through the Lyapunov coefficients, which are obtained by the expansion of the Poincaré first return map in power series, the details of which are given below.

The quasi-cylindrical coordinates $u=\rho \cos\theta, v=\rho \sin\theta$, and $w=\rho\omega$ was successfully used in \cite{Aulb85} to solve the center-focus problem and later by \cite{Bui14} to obtain a reduced form of the Poincaré map and also to study the cyclicity of Hopf points \cite{Gar18}. Rewriting system \eqref{formnorm} in these coordinates, we have
 \begin{equation*}
	\begin{array}{l}
		\dot{\rho} = P(\rho\cos\theta,\rho\sin\theta,\rho\omega)\cos\theta+Q(\rho\cos\theta,\rho\sin\theta,\rho\omega)\sin\theta ,\\
		\dot{\theta} = 1+\rho^{-1}\left(Q(\rho\cos\theta,\rho\sin\theta,\rho\omega)\cos\theta-P(\rho\cos\theta,\rho\sin\theta,\rho\omega)\sin\theta\right), \\
		\dot{\omega} = \lambda \omega + \rho^{-1} \left(R(\rho\cos\theta,\rho\sin\theta,\rho\omega)-\omega\dot{\rho}\right).
	\end{array}
\end{equation*}
Since $P$, $Q$, and $R$ are analytic functions with Taylor series expansions containing terms of degree two or higher, we can introduce the notation $\tilde{A}(\rho,\theta,\omega)=\rho^{-2} A(\rho \cos\theta,\rho \sin\theta,\rho\omega)$ for $A=P,Q,R$. Thus, noting that the term $\dot{\rho}$ in the equation for $\dot{\omega}$ can be rewritten as $\dot{\rho}=\rho^2\left(\tilde{P}\cos\theta+\tilde{Q}\sin\theta\right)$, the system above becomes
\begin{equation}
	\label{formpolar}
	\begin{array}{l}
		\dot{\rho}= \rho^2\left(\tilde{P}(\rho,\theta,\omega)\cos\theta+\tilde{Q}(\rho,\theta,\omega)\sin\theta\right),\\
		\dot{\theta} = 1+\rho\left(\tilde{Q}(\rho,\theta,\omega)\cos\theta-\tilde{P}(\rho,\theta,\omega)\sin\theta\right), \\
		\dot{\omega} = \lambda \omega + \rho \left(\tilde{R}(\rho,\theta,\omega)-\omega\left(\tilde{P}(\rho,\theta,\omega)\cos\theta+\tilde{Q}(\rho,\theta,\omega)\sin\theta\right)\right).
	\end{array}
\end{equation} 
Now, taking $\theta$ as the new independent variable, we obtain
\begin{equation}
	\label{dr/do}
	\begin{array}{l}
		\dfrac{d\rho}{d\theta} = \dfrac{\rho^2\left(\tilde{P}\cos\theta+\tilde{Q}\sin\theta\right)}{1+\rho(\tilde{Q}\cos\theta-\tilde{P}\sin\theta)} = \mathcal{R}(\theta,\rho,\omega) ,\\
		\dfrac{d\omega}{d\theta} = \dfrac{\lambda \omega +\rho(\tilde{R}-\omega(\tilde{P}\cos\theta+\tilde{Q}\sin\theta))}{1+\rho(\tilde{Q}\cos\theta-\tilde{P}\sin\theta)} =\lambda \omega+ \Omega(\theta,\rho,\omega),
	\end{array}
\end{equation}
 with $\mathcal{R}$ and $\Omega$ analytical periodic functions in the cylinder $C=\{ (\rho,\omega) : |\rho|\leq r^*,\; |\omega|\leq h^* \},$ for $r^*$ and $h^*$ small enough. For $(\rho_0,\omega_0)\in C$, denote by $\varphi(\theta,\rho_0,\omega_0)=\left(\rho(\theta,\rho_0,\omega_0),\omega(\theta,\rho_0,\omega_0)\right)$ the solution of \eqref{dr/do} with the initial condition $\left(\rho(0,\rho_0,\omega_0),\omega(0,\rho_0,\omega_0)\right) = (\rho_0,\omega_0)$, which is also analytic for $0\leq\theta\leq2\pi$ and $(\rho_0,\omega_0)\in C$. Furthermore, $\rho(\theta,0,0)=\omega(\theta,0,0)\equiv 0$, since the origin is a singular point of \eqref{dr/do}. The next definitions and lemma can be seen in \cite{Bui14}.
\begin{definition}
	We define the {\it Poincaré map} associated with the system \eqref{formnorm} by
	\begin{equation*}
		\Pi(\rho_0,\omega_0)=\varphi(2\pi,\rho_0,\omega_0).
	\end{equation*}
\end{definition}

\begin{definition}
	For the system \eqref{formnorm}, the map defined by
	\begin{eqnarray*}
		d(\rho_0,\omega_0) & = & \Pi(\rho_0,\omega_0) - (\rho_0,\omega_0) = \varphi(2\pi,\rho_0,\omega_0)-(\rho_0,\omega_0) \\
		& = & \left(\rho(2\pi,\rho_0,\omega_0)-\rho_0,\omega(2\pi,\rho_0,\omega_0)-\omega_0\right) = \left(d_1(\rho_0,\omega_0),d_2(\rho_0,\omega_0)\right)
	\end{eqnarray*}
	is called the {\it displacement function}.
\end{definition}
\begin{lemma}\label{lemad2}
	Let $d(\rho_0,\omega_0)$ be the displacement function associated with the system \eqref{formnorm}. Then, there is a unique analytical function $\tilde{\omega}$, defined in a neighborhood $V$ of $\rho_0=0$, such that $d_2(\rho_0,\tilde{\omega}(\rho_0))=0$.
\end{lemma}
Thus, we define the {\it reduced displacement function}, given by
\begin{equation}
	\label{valfocal}
	\overline{d}(\rho_0)=d_1(\rho_0,\tilde{\omega}(\rho_0))=\sum_{j\geq 1} v_j \rho_0^j.
\end{equation}
We have that the first nonzero coefficient $v_j$ from \eqref{valfocal} correspond to an odd power in $\rho_0$, see \cite[Theorem 2]{Bui14} . Therefore, we can state the following definition.
\begin{definition}
	For system \eqref{formnorm} and its values $v_j$ given in \eqref{valfocal}, the coefficients
	\begin{equation}
		\label{eqLyap}
		\ell_k=v_{2k+1},
	\end{equation}
	where $k\geq 0$, are called Lyapunov coefficients.
\end{definition}
Note that, the isolated zeros of the reduced displacement function correspond to the limit cycles of system \eqref{formnorm}, and when this function is identically null, we have a center on the center manifold. That is, if $\ell_k=0$ for all $k$, then the origin of \eqref{formnorm} is a center on the center manifold, and the converse is also true.

Theorem 7, found in \cite{Gar18}, relates the Lyapunov coefficients $\ell_k$, given in \eqref{eqLyap}, to the focus quantities $L_k$ in \eqref{eqFocal}. More precisely, given a positive integer $k > 2$, it holds that
\begin{equation*}
	\ell_1 = \pi L_1 \quad \text{and} \quad \ell_{k-1} = \pi L_{k-1} \quad\textrm{mod}\left\langle L_1, L_2, \dots, L_{k-2}\right\rangle.
\end{equation*}
Thus, we can consider one or the other in the study of the center-focus problem and the cyclicity from Hopf points. But, for computational purposes, we give preference to the focus quantities.

\section{Isochronicity Problem}\label{secIsoc}
Consider $\varphi(\te,\rz,\wz)=(\rho(\te,\rz,\wz),\w(\te,\rz,\wz))$ to be the solution of system \eqref{dr/do}. Its Taylor series expansion around $\rz=0$ is given by
\begin{equation}
	\label{sol}
	\begin{array}{l}
		\rho(\te,\rz,\wz) = \sum_{i\geq 0} u_{i}(\te,\wz) \rz^{i+1} = u_0(\te,\wz)\rz +u_1(\te,\wz)\rz^2+\cdots \\ \\
		\w(\te,\rz,\wz) = \sum_{i\geq 0} v_{i}(\te,\wz) \rz^{i} = v_0(\te,\wz) +v_1(\te,\wz)\rz+\cdots
	\end{array}
\end{equation}
where $u_{0}(0,\wz) = 1$, $v_0(0,\wz)=\wz$, and $u_i(0,\wz)=v_i(0,\wz)=0$ for $i\geq1$. The functions $u_i$'s and $v_i$'s are uniquely determined since \eqref{sol} is the solution of a Cauchy problem. According to Lemma \ref{lemad2}, there exists a unique analytic function $\tilde{\w}(\rz)$ in a neighborhood of $\rz=0$ such that $d_2(\rz,\tilde{\w}(\rz))\equiv0$. Now, let us consider the expression for $\dot{\te}$. Substituting the expansion of solution $\varphi(\te,\rz,\tilde{\w}(\rz))=\left(\rho(\te,\rz,\tilde{\w}(\rz)),\w(\te,\rz,\tilde{\w}(\rz))\right)$ given in \eqref{sol} into the expression for $\dot{\te}$ in \eqref{formpolar} and performing a Taylor series expansion around $\rz=0$, we obtain
\begin{equation}
	\label{dt/do}
	\dfrac{dt}{d\te}= \dfrac{1}{1+\sum_{k=1}^{\infty}F_k(\te)\rz^k}=1+\sum_{k=1}^{\infty}\Psi_k(\te) \rz^k.
\end{equation}
Integrating \eqref{dt/do}, since $\sum_{k=1}^{\infty}\Psi_k(\te) \rz^k$ is analytic, we obtain
\begin{equation}
	t=\int_{0}^{\te}\left(1+\sum_{k=1}^{\infty}\Psi_k(\tau) \rz^k\right)d\tau=\te+\sum_{k=1}^{\infty}\varPhi_k(\te)\rz^k,
\end{equation}
where $\varPhi_k(\te)= \int_{0}^{\te} \Psi_k(\tau)d\tau$.
Since $\rz$ is fixed and small enough, the above series is convergent for $\te\in[0,2\pi]$. Hence, taking $\te=2\pi$, we obtain the period function for a system in $\mathbb{R}^3$ without the initial need to restrict to a center manifold $W^c$, which is given by
\begin{equation}
	\label{funcperiod}
	T(\rz)=2\pi\left(1+\sum_{k=1}^{\infty}T_{2k}\rz^{2k}\right),
\end{equation}
with $T_{2k}=\frac{1}{2\pi}\varPhi_{2k}(2\pi)$, $k\geq1$, being called the \textit{isochronicity constants}.

A center on the center manifold of system \eqref{formnorm} is said to be isochronous if and only if its period function is constant. That is, if all its coefficients $T_{2k}$ are null for $k\geq1$.

\section{Proofs of Theorems \ref{Teo1} and \ref{Teo2}}\label{secProof}
	
\subsection{Proof of Theorem \ref{Teo1}}

 The Jacobian matrix of system \eqref{sist3Dchaotic} evaluated at its singular point $E_1$ is given by
\begin{equation*}
	J(E_1) = \begin{pmatrix}
		-a & a+\dfrac{1}{d} & 0 \\
		b-\dfrac{1}{d} & c & 0 \\
		0 & 0 & -d
	\end{pmatrix},
\end{equation*}
whose characteristic polynomial is $P(\lambda)=\lambda^3+\alpha\lambda^2+\beta \lambda+\gamma$, with
\begin{eqnarray*}
	\alpha&=&a-c+d;\\
	\beta&=&-\dfrac{-1 - a d + b d + a b d^2 + a c d^2 - a d^3 + c d^3}{d^2};\\
	\gamma&=&-\dfrac{-1 - a d + b d + a b d^2 + a c d^2}{d}.
\end{eqnarray*}

Since we are looking for centers on the center manifold that are also Hopf points, the Jacobian matrix $J(E_1)$ must have a pair of purely imaginary eigenvalues and a nonzero real eigenvalue. It follows from the characteristic polynomial $P(\lambda)$ that the singular point $E_1$ is a Hopf point if and only if  $\beta>0$ and $\gamma-\alpha\beta=0$. This leads us to the conditions $a=c$ and $(1+cd)(1-bd)- c^2 d^2>0$ and, in this case, the eigenvalues are as follows $\lambda_{1,2}=\pm \dfrac{\sqrt{(1+cd)(1-bd)- c^2 d^2}}{d}i$ and $\lambda_3=-d$, proving the first statement of the theorem. 
%Now, translating the singular point $E_1$ to the origin and introducing a new parameter $k>0$ given by the equation
%\begin{equation*}
%	1 - b d + c d - b c d^2 - c^2 d^2-k^2=0,
%\end{equation*}
%that is,
%\begin{equation}\label{eqbE1}
%	b=\dfrac{1 + c d - c^2 d^2 - k^2}{d(1+c d)},
%\end{equation}
%in order to simplify the expressions of the eigenvalues and eigenvectors, we obtain the system

Now, introducing a new real parameter $k>0$ given by the equation
	\begin{equation*}
		(1+cd)(1-bd)- c^2 d^2 -k^2=0,
	\end{equation*}
	that is,
	\begin{equation}\label{eqb}
		b=\dfrac{1 + c d - c^2 d^2 - k^2}{d(1+c d)},
	\end{equation}
	we can simplify the expressions of the eigenvalues and eigenvectors. This expression for $b$ is well defined, as  $d\neq0$ by hypothesis, and the condition for the singular point $E_1$ to be a Hopf point implies that $1+cd\neq0$. Considering this and translating the singular point $E_1$ to the origin, we can write system \eqref{sist3Dchaotic} as
\begin{equation}\label{sistE1}
	\begin{array}{l}
		\dot{x}=-c x + c y + \dfrac{y}{d} + y z,\\ \\
		\dot{y}=-\dfrac{c^2 d^2 x}{c d^2+d}-\dfrac{k^2 x}{c d^2+d}+\dfrac{c d x}{c d^2+d}+\dfrac{x}{c d^2+d}+c y-\dfrac{x}{d}-x z,\\ \\
		\dot{z}=x y - d z.
	\end{array}
\end{equation}

The eigenvalues of linear part of system \eqref{sistE1} are $\lambda_{1,2}=\pm\dfrac{k}{d}i$ and $\lambda_3=-d$ with the eigenvectors given by $v_{1,2}=\left(\dfrac{1+c d}{c d\pm i k},1,0\right)$ and $v_3=(0,0,1)$. Thus, introducing the change of coordinates
\begin{equation*}
	(x,y,z)\mapsto\left(\dfrac{c d k+k}{c^2 d^2+k^2}u+\dfrac{c d (c d+1)}{c^2 d^2+k^2}v,v,w\right),
\end{equation*}
and a time rescaling $\tau=\frac{k}{d}t$, system \eqref{sistE1} becomes
\begin{equation}\label{sistE1norm}
	\begin{array}{l}
		\dot{u}=v+\dfrac{c d^2 (c d+1)}{c^2 d^2 k+k^3}u w+\dfrac{d \left(2 c^4 d^4+2 c^3 d^3+2 c^2 d^2 k^2+c^2 d^2+k^4\right)}{k^2 (c d+1) \left(c^2 d^2+k^2\right)}v w,\\ \\
		\dot{v}=-u-\dfrac{d (c d+1)}{c^2 d^2+k^2}u w-\dfrac{c d^2 (c d+1)}{k \left(c^2 d^2+k^2\right)}v w,\\ \\
		\dot{w}=-\dfrac{d^2}{k}w+\dfrac{d (c d+1)}{c^2 d^2+k^2}u v+\dfrac{c d^2 (c d+1)}{k \left(c^2 d^2+k^2\right)}v^2.
	\end{array}
\end{equation}

Applying the algorithm of Section \ref{secCenter}, we can compute the focus quantities $L_k$, $k=1,2,3$, which are
\allowdisplaybreaks
\begin{align*}
	L_1&=d \left(k^2+4 c^2-1\right)+2 \left(k^2+1\right) c+2 c \left(2 c^2-1\right) d^2;\\[1ex]
	L_2&=c \big(12 k^6 (2 c+d)+k^4 \left(64 c^3 d^2+8 c^2 \left(3 d^2+4\right) d+2 c \left(d^4+8\right)+d^5\right)\\
	&\quad+2 k^2 \left(48 c^5 d^4+96 c^4 d^3-24 c^3 d^2 \left(d^2-4\right)+c^2 d \left(d^6-4 d^4-36 d^2+48\right)-2 c \left(d^4+12 d^2-6\right)-6 d\right)\\
	&\quad+d^4 \left(8 c^5 d^4+16 c^4 d^3-4 c^3 d^2 \left(d^2-4\right)+c^2 \left(8 d-6 d^3\right)+c \left(2-4 d^2\right)-d\right)\big),
\end{align*}
and $L_3$ can be seen in Appendix \ref{apendice}. Using the routine min\textit{AssGTZ} \cite{minAssGTZRoutine} of \textsc{Singular} \cite{Singular}, we obtain the decomposition of the radical of the ideal generated by the
three focus quantities $\langle L_1,L_2,L_3\rangle$ into an intersection of the following minimal prime ideals: 
\begin{align*}
		I_1&=\langle 4 c^3 d^2 + 4 c^2 d - 2 cd^2 + 2 c - d, k \rangle;\\
		I_2&=\big\langle d^4 + 4 k^2, -d^3 k - d^3 + 8 c k^2 + 2 d k^2 + 2 d k, 4 c d k + d^2 k + d^2 + 2 k^2 + 2 k, c d^2 + 2 c k + d k + d, \\
		&\qquad 4 c^2 k - k^2 - 2 k - 1 \rangle;\\
		I_3&=\langle d^4 + 4 k^2, d^3 k - d^3 + 8 c k^2 + 2 d k^2 - 2 d k, 4 c d k + d^2 k - d^2 - 2 k^2 + 2 k, c d^2 - 2 c k - d k + d, \\
		&\qquad4 c^2 k + k^2 - 2 k + 1 \rangle;\\
		I_4&=\langle d^2 - 2, 2 c + d \rangle;\\
		I_5&=\langle k + 1, c \rangle;\\
		I_6&= \langle k-1, c \rangle;\\
		I_7&=\langle d, c \rangle;\\
		I_8&=\langle k^2+1, 1+c d \rangle;\\
		I_9&=\langle d^4+4k^2, 2c-d\rangle.\\
\end{align*}
From these ideals, we can find necessary conditions for the origin to be a center on the center manifold, which are obtained vanishing the generators of $I_j$, for $j=1,2,\dots,9$. Note that $k$ was assumed to be positive, while $d$ is nonzero by the hypothesis. Therefore, the generators of the ideals $I_1, I_2, I_3, I_5, I_7, I_8$ and $I_9$ do not vanish and we can discard them. From $I_4\equiv0$, we get $d=\sqrt{2}$ and $c=-1/\sqrt{2}$, which implies that $1+cd=0$. However, this contradicts the condition $1+cd\neq0$ noted after equation \eqref{eqb}. Finally, from $I_6\equiv0$, we have $k=1$ and $c=0$. Substituting these conditions into system \eqref{sistE1norm}, the system becomes
\begin{equation}\label{sistE1cond6}
	\begin{array}{l}
		\dot{u}=v+d v w,\\
		\dot{v}=-u-d u w, \\
		\dot{w}=-d^2 w+d u v.
	\end{array}
\end{equation}
Note that $H(u,v,w)=u^2+v^2$ is a first integral for this system, because $\langle (\dot{u},\dot{v},\dot{w}), \nabla H\rangle=0$. Thus, the origin of system \eqref{sistE1norm} is a center on the center manifold if and only if $k=1$ and $c=0$. Substituting these conditions back into the original system's parameters using equation \eqref{eqb}, we obtain the center conditions $b=c=0$.

Now, we will prove that system \eqref{sistE1cond6} does not have an isochronous center at the origin on the center manifold. Proceeding as in Section \ref{secIsoc}, we can obtain the expansion of the period function, given by
	\begin{equation*}
		T(\rho_0)=2\pi \left(-1-\dfrac{d^4}{8 \left(d^4+4\right)}\rho_0^4+\cdots\right).
	\end{equation*}
	Therefore, since the constant $T_4=-\dfrac{d^4}{8 \left(d^4+4\right)}$ only vanishes for $d=0$, and $d$ is nonzero by hypothesis, the center is never isochronous. \qed

\begin{remark}
For points $E_2^\pm$ and $E_3^\pm$, the computations required to obtain the normal forms of the system, after translating the points to the origin, are very hard. In fact, we are not even able to determine the general conditions under which they would be a Hopf points. Therefore, we do not study these points in their full generality. Instead, for specific values of $a$ and $b$, we present conditions on the parameters $c$ and $d$ under which they are Hopf points. Even under these conditions, however, analyzing the center-focus problem remains unfeasible. Nevertheless, we provide a specific choice of parameters $c$ and $d$ for which the points $E_2\pm$ and $E_3\pm$ are foci on the center manifolds. The results are presented below.
\end{remark}

\begin{propo}\label{propo}
	Consider system \eqref{sist3Dchaotic}, with $d\neq0$. Let $a=1$ and $b=-1$ (resp. $a=-1$ and $b=1$). If $c>1$ and $d = \frac{1}{2 (-1 + c)} \left( 2 + c + 4 c^{3/2} + c^2 - \sqrt{c^{3/2} \left(16 + 9 \sqrt{c} + 8 c + 22 c^{3/2} + 8 c^2 + c^{5/2}\right)} \right)$ \Bigg(resp. $c<-1$ and $d = \frac{1}{2 (-1 - c)} \left( -2 + c + 4 c \sqrt{-c} - c^2 - \sqrt{-4 (-1 - c) (1 - 2 c + c^2) + \left(2 - c - 4 c \sqrt{-c} + c^2\right)^2} \right)$\Bigg), then the singular points $E_2^\pm$ (resp. $E_3^\pm$) are a pair of Hopf points of the system.
\end{propo}
	
	\begin{proof}
		Consider system \eqref{sist3Dchaotic} with $d\neq0$, $a=1$ and $b=-1$, and the points $E_2^\pm$. Since the singular points are symmetric, we will consider only $E_2^-$. The characteristic polynomial associated to the Jacobian matrix of the system evaluated at $E_2^-$ is $P(\lambda)=\lambda^3+\alpha \lambda^2+\beta \lambda + \gamma$, with
		\begin{align*}
			\alpha &=1-c+d ,\\
			\beta &= -\frac{1}{\sqrt{c}} + \sqrt{c} - \frac{d}{\sqrt{c}} + \sqrt{c} d , \\
			\gamma &= 4 \sqrt{c} \left(1 + d + \sqrt{c} d\right).
		\end{align*}

Thus, $E_2^-$ is a Hopf point if $\beta > 0$ and $\gamma-\alpha\beta = 0$. Using symbolic computation software such as Mathematica, we obtain the following conditions: $c > 1$ and $$d = \frac{1}{2 (-1 + c)} \left( 2 + c + 4 c^{3/2} + c^2 - \sqrt{c^{3/2} \left(16 + 9 \sqrt{c} + 8 c + 22 c^{3/2} + 8 c^2 + c^{5/2}\right)} \right).$$

Now, consider the points $E_3^\pm$ with $a = -1$ and $b = 1$. Again, we will only consider the point $E_3^-$. The characteristic polynomial associated is $P(\lambda)=\lambda^3+\alpha \lambda^2+\beta \lambda + \gamma$, with
\begin{align*}
	\alpha &=-1-c+d ,\\
	\beta &= \sqrt{-c} + \frac{\sqrt{-c}}{c} - d\sqrt{-c}  - d \frac{\sqrt{-c}}{c} , \\
	\gamma &=4 \left( \sqrt{-c} + \left( -\sqrt{-c} + c \right) d \right).
\end{align*}
Thus, we obtain as a condition for the point to be a Hopf point, $c<-1$ and
$$d = \frac{1}{2 (-1 - c)} \left( -2 + c + 4 c \sqrt{-c} - c^2 - \sqrt{-4 (-1 - c) (1 - 2 c + c^2) + \left(2 - c - 4 c \sqrt{-c} + c^2\right)^2} \right).$$
	\end{proof}
	
\begin{corol}\label{coro}
	Under the conditions of Proposition \ref{propo}, and for $c=9$ and $d=\frac{5}{2}(5+3\sqrt{3})$ \big(resp. $c=-9$ and $d=-\frac{5}{2}(5+3\sqrt{3})$\big), the points $E_2^\pm$ (resp. $E_3^\pm$) are foci on their center manifolds.
\end{corol}

\begin{proof}
	Let us consider the points $E^\pm_2$. Under the conditions of Proposition \ref{propo}, we then have $a=1$ and $b=-1$. For $c=9$, the Hopf condition yields $d=\frac{5}{2}(5+3\sqrt{3})$. Due to the symmetry of the points, we will consider the point $E_2^-$. By translating this point to the origin and applying a change of coordinates combined with a time rescaling, the system can be written as
\begin{equation}\label{sistE2norm}
	\begin{aligned}
		u' &= -v - \frac{399 + 248 \sqrt{3}}{12980} u^2 - \frac{4 \sqrt{4491 - 2425 \sqrt{3}}}{649} u v - \frac{\sqrt{14277297 - 7644019 \sqrt{3}}}{9735} u w \\
		&\quad + \frac{1011 - 167 \sqrt{3}}{5192} v^2 + \frac{-14653 + 8809 \sqrt{3}}{7788} v w + \frac{2 \left(-9901 + 5663 \sqrt{3}\right)}{9735} w^2, \\[1ex]
		v' &= u - \frac{\sqrt{430587 - 246009 \sqrt{3}}}{3245} u^2 + \frac{15 \left(-29 + 8 \sqrt{3}\right)}{1298} u v + \frac{11881 - \frac{22544}{\sqrt{3}}}{6490} u w \\
		&\quad - \frac{\sqrt{\frac{1}{2} \left(1707 - 359 \sqrt{3}\right)}}{649} v^2 - \frac{\sqrt{2228043 - 1245841 \sqrt{3}}}{1947} v w + \frac{9408 - 7474 \sqrt{3}}{9735 \sqrt{9 + 5 \sqrt{3}}} w^2, \\[1ex]
		w' &= \frac{3}{4} \sqrt{201 - 115 \sqrt{3}} w + \frac{3 \sqrt{3 \left(7401 + 4133 \sqrt{3}\right)}}{4720} u^2 + \frac{63 + 5 \sqrt{3}}{236} u v \\
		&\quad + \frac{117 - 16 \sqrt{3}}{295} u w - \frac{\sqrt{\frac{1}{2} \left(1881 - 53 \sqrt{3}\right)}}{944} v^2 + \frac{\sqrt{438 - 190 \sqrt{3}}}{59} v w + \frac{3 \left(48 + \sqrt{3}\right)}{295 \sqrt{9 + 5 \sqrt{3}}} w^2.
	\end{aligned}
\end{equation}
	Applying the algorithm of Section \ref{secCenter}, we obtain
	\begin{equation*}
		L_1=\frac{2652236723}{40 \sqrt{3 \left(79295001361791115851 + 40605175853700180095 \sqrt{3}\right)}} \neq 0.
	\end{equation*}
	Hence, the points $E_2^\pm$ are foci on their center manifolds.
	
	Now, consider the points $E_3^\pm$ under the conditions of Proposition \ref{propo}, with $a=-1$ and $b=1$. For $c=-9$, we obtain $d=-\frac{5}{2}(5+3\sqrt{3})$. Again, we will consider only the point $E_3^-$. After a translation, a change of coordinates, and a time rescaling, we can write the system in the following form
	\begin{equation}\label{sistE3norm}
		\begin{aligned}
			u' &= -v + \frac{399 + 248 \sqrt{3}}{12980} u^2 - \frac{4 \sqrt{4491 - 2425 \sqrt{3}}}{649} u v - \frac{\sqrt{14277297 - 7644019 \sqrt{3}}}{9735} u w \\
			&\quad + \frac{-1011 + 167 \sqrt{3}}{5192} v^2 + \frac{14653 - 8809 \sqrt{3}}{7788} v w + \frac{2 \left(9901 - 5663 \sqrt{3}\right)}{9735} w^2, \\[1ex]
			v' &= u - \frac{\sqrt{430587 - 246009 \sqrt{3}}}{3245} u^2 + \frac{15 \left(29 - 8 \sqrt{3}\right)}{1298} u v + \frac{-35643 + 22544 \sqrt{3}}{19470} u w \\
			&\quad - \frac{\sqrt{\frac{1}{2} \left(1707 - 359 \sqrt{3}\right)}}{649} v^2 - \frac{\sqrt{2228043 - 1245841 \sqrt{3}}}{1947} v w + \frac{9408 - 7474 \sqrt{3}}{9735 \sqrt{9 + 5 \sqrt{3}}} w^2, \\[1ex]
			w' &= \frac{3}{4} \sqrt{51 - 25 \sqrt{3}} w + \frac{3 \sqrt{3 \left(7401 + 4133 \sqrt{3}\right)}}{4720} u^2 + \frac{-63 - 5 \sqrt{3}}{236} u v \\
			&\quad + \frac{-117 + 16 \sqrt{3}}{295} u w - \frac{\sqrt{\frac{1}{2} \left(1881 - 53 \sqrt{3}\right)}}{944} v^2 + \frac{\sqrt{438 - 190 \sqrt{3}}}{59} v w + \frac{3 \left(48 + \sqrt{3}\right)}{295 \sqrt{9 + 5 \sqrt{3}}} w^2.
		\end{aligned}
	\end{equation}
	
	In this case, we obtain
	\begin{equation*}
		L_1=\frac{\sqrt{44674468353 - 21555759691 \sqrt{3}}}{38940 \left(523 - 225 \sqrt{3}\right)} \neq 0.
	\end{equation*}
	Thus, the points $E_3^\pm$ are foci on their center manifolds.
\end{proof}

\subsection{Proof of Theorem \ref{Teo2}}

 Let us first consider the pair $E_4^\pm$. Since it is a pair of symmetric points, we will analyze only one of them, i.e., we consider point $E_4^-$.

By expression \eqref{ptsE4E5}, the characteristic polynomial of Jacobian matrix $J(E_4^-)$ is  $P(\lambda)=\lambda^3+\alpha\lambda^2+\beta \lambda+\gamma$, where
\begin{eqnarray*}
	\alpha&=&a-c;\\
	\beta&=&-\dfrac{a^2+a b +a c + b c-a\sqrt{(a+b)^2+4a c}+c\sqrt{(a+b)^2+4a c}}{2a c};\\
	\gamma&=&-2\sqrt{(a+b)^2+4a c}.
\end{eqnarray*}
Imposing the conditions for the existence of a Hopf point, $\beta>0$ and $\gamma-\alpha\beta=0$, we have that $a=-c$, $c>0$ and $b>3c$. Hence, the eigenvalues are $\lambda_{1,2}=\pm\dfrac{\sqrt[4]{(b-3c)(b+c)}}{\sqrt{c}}i$ and $\lambda_3=2c$.
By translating point $E_4^-$ to the origin and inserting a new real parameter $h>0$ by the equation
\begin{equation*}
	\sqrt{(b-3 c) (b+c)}+b-c-h^2=0,
\end{equation*}
so that
\begin{equation*}
	b=\dfrac{4 c^2+2 c h^2+h^4}{2 h^2},
\end{equation*}
we can write system \eqref{sist3Dchaotic} as
\begin{equation}\label{sistE4m}
	\begin{array}{l}
		\dot{x}=\dfrac{\sqrt{2} \sqrt{c} }{h}z+c x+\dfrac{h^2}{2}y+y z,\\ \\
		\dot{y}=\dfrac{2 c^2}{h^2}x+\dfrac{h}{\sqrt{2} \sqrt{c}}z+c y-x z, \\ \\
		\dot{z}=\dfrac{\sqrt{2} \sqrt{c}}{h}x-\dfrac{h}{\sqrt{2} \sqrt{c}}y+x y,
	\end{array}
\end{equation}
with conditions $c>0$, $h>0$ and $h^4-4c^2>0$. The linear part of system \eqref{sistE4m} has
eigenvalues $\lambda_{1,2}=\pm\dfrac{\sqrt{h^4-4 c^2}}{\sqrt{2c} h}i$ and $\lambda_3=2c$ and eigenvectors $v_1=\overline{v_2}=\left(\dfrac{-h\sqrt{c} \sqrt{h^4-4 c^2}-2 \sqrt{2} c\;i}{\sqrt{2} \sqrt{h^4-4 c^2}+4 c^{3/2} h\;i},\dfrac{2 c^{3/2} \sqrt{h^4-4 c^2}- \sqrt{2} h^3\;i}{\sqrt{2} h \sqrt{h^4-4 c^2}+4 c^{3/2} h^2\;i},1\right)$ and $v_3=\left(\dfrac{h^2}{2c},1,0\right)$. Now, performing the change of coordinates
\begin{equation*}
	\begin{aligned}
		(x,y,z)\mapsto&\Bigg(-\dfrac{2 c \left(c h^2-1\right) \sqrt{h^4-4 c^2}}{8 c^3 h^2-4 c^2+h^4}u-\dfrac{ \sqrt{c} h \left(4 c^2+h^4\right)}{\sqrt{2} \left(8 c^3 h^2-4 c^2+h^4\right)}v+\dfrac{h^2}{2 c}w,\\
		&\quad\dfrac{ \left(4 c^3+h^2\right) \sqrt{h^4-4 c^2}}{8 c^3 h^2-4 c^2+h^4}u-\dfrac{\sqrt{2} c^{3/2} \left(4 c^2+h^4\right)}{8 c^3 h^3-4 c^2 h+h^5}v+w, v\Bigg)
	\end{aligned}
\end{equation*}
and a time rescaling $\tau=\dfrac{\sqrt{h^4-4 c^2} }{\sqrt{2} \sqrt{c} h}t$, system \eqref{sistE4m} becomes
\begin{equation}\label{sistE4mnorm}
	\begin{aligned}
		\dot{u}={}&v+\dfrac{h \left(4 c^2+h^4\right)}{\sqrt{2} \sqrt{c} \left(4 c^2-h^4\right)}v w-\dfrac{c \left(4 c^2+h^4\right)^2}{\left(4 c^2-h^4\right) \left(8 c^3 h^2-4 c^2+h^4\right)}v^2+\dfrac{2 \sqrt{2} c^{3/2} h \left(-4 c^3+c h^4-2 h^2\right)}{\sqrt{h^4-4 c^2} \left(8 c^3 h^2-4 c^2+h^4\right)}u v,\\ \\
		\dot{v}={}& -u+\dfrac{h^3 }{\sqrt{2} \sqrt{c} \sqrt{h^4-4 c^2}}w^2-\dfrac{c \left(4 c^2+h^4\right)^2}{\left(8 c^3 h^2-4 c^2+h^4\right)^2}u v+\dfrac{h \left(4 c^2+h^4\right)}{\sqrt{2} \sqrt{c} \left(8 c^3 h^2-4 c^2+h^4\right)}u w\\
		& {}-\dfrac{2 c \left(4 c^2 h^2+h^6\right)}{\sqrt{h^4-4 c^2} \left(8 c^3 h^2-4 c^2+h^4\right)}v w+\dfrac{2 \sqrt{2} c^{3/2} h \left(2 c-h^2\right) \left(2 c+h^2\right) \left(4 c^3+h^2\right) \left(c h^2-1\right)}{\sqrt{h^4-4 c^2} \left(8 c^3 h^2-4 c^2+h^4\right)^2}u^2\\
		& {}+\dfrac{\sqrt{2} c^{5/2} h \left(4 c^2+h^4\right)^2}{\sqrt{h^4-4 c^2} \left(8 c^3 h^2-4 c^2+h^4\right)^2}v^2, \\ \\
		\dot{w}={}&\dfrac{2 \sqrt{2} c^{3/2} h }{\sqrt{h^4-4 c^2}}w+\dfrac{4 c^3 \left(2 c-h^2\right) \left(2 c+h^2\right) \left(4 c^3+h^2\right) \left(c h^2-1\right) \left(4 c^2+h^4\right)}{\sqrt{h^4-4 c^2} \left(8 c^3 h^2-4 c^2+h^4\right)^3}u^2 \\
		&{}+\dfrac{2 c^3  \left(4 c^2+h^4\right) \left(64 c^6 h^2-48 c^5-16 c^4 h^6+40 c^3 h^4-16 c^2 h^2-3 c h^8+4 h^6\right)}{\sqrt{h^4-4 c^2} \left(8 c^3 h^2-4 c^2+h^4\right)^3}v^2\\
		&{}+\dfrac{c h^2 \left(4 c^2+h^4\right)}{\sqrt{h^4-4 c^2} \left(8 c^3 h^2-4 c^2+h^4\right)} w^2+\dfrac{4 \sqrt{2} c^{3/2} h \left(2 c-h^2\right) \left(2 c+h^2\right) \left(4 c^3+h^2\right) \left(c h^2-1\right)}{\sqrt{h^4-4 c^2} \left(8 c^3 h^2-4 c^2+h^4\right)^2} v w\\
		&{}+\dfrac{\sqrt{2} c^{3/2} \left(4 c^2 + h^4\right) \left(-16 c^5 - 8 c^2 (1 + 4 c^4) h^2 + 8 c^3 (1 + 8 c^4) h^4 + 2 (1 + 4 c^4) h^6 - c h^8\right)}{h \left(8 c^3 h^2-4 c^2+h^4\right)^3} u v\\
		&{}+\frac{c  \left(4 c^2+h^4\right)^2}{\left(8 c^3 h^2-4 c^2+h^4\right)^2}u w,
	\end{aligned}
\end{equation}
Using the algorithm of Section \ref{secCenter} for obtaining the focus quantities, we have
\begin{eqnarray*}
	L_1&=&-h c^{7/2} \sqrt{h^4-4 c^2} \left(h^4+4 c^2\right)^2.
%	L_2&=&-27 h^{35} c^{9/2}+264 h^{35} c^{17/2}+579 h^{33} c^{15/2}+233 h^{33} c^{23/2}-966 h^{31} c^{13/2}+11092 h^{31} c^{21/2}\\
%	&&-400 h^{31} c^{29/2}-18304 h^{29} c^{19/2}+37012 h^{29} c^{27/2}+320 h^{29} c^{35/2}+9456 h^{27} c^{17/2}\\
%	&&-76512 h^{27} c^{25/2}+18816 h^{27} c^{33/2}+54688 h^{25} c^{23/2}+30096 h^{25} c^{31/2}-1280 h^{25} c^{39/2}\\
%	&&+15456 h^{23} c^{21/2}-177472 h^{23} c^{29/2}+129280 h^{23} c^{37/2}+548672 h^{21} c^{27/2}-456896 h^{21} c^{35/2}\\
%	&&-10240 h^{21} c^{43/2}-288768 h^{19} c^{25/2}+2313216 h^{19} c^{33/2}+380928 h^{19} c^{41/2}-2194688 h^{17} c^{31/2}\\
%	&&+1827584 h^{17} c^{39/2}+40960 h^{17} c^{47/2}+247296 h^{15} c^{29/2}-2839552 h^{15} c^{37/2}+2068480 h^{15} c^{45/2}\\
%	&&-3500032 h^{13} c^{35/2}-1926144 h^{13} c^{43/2}+81920 h^{13} c^{51/2}+2420736 h^{11} c^{33/2}-19587072 h^{11} c^{41/2}\\
%	&&+4816896 h^{11} c^{49/2}+18743296 h^9 c^{39/2}-37900288 h^9 c^{47/2}-327680 h^9 c^{55/2}-3956736 h^7 c^{37/2}\\
%	&&+45432832 h^7 c^{45/2}-1638400 h^7 c^{53/2}-9486336 h^5 c^{43/2}-3817472 h^5 c^{51/2}-1769472 h^3 c^{41/2}\\
%	&&+17301504 h^3 c^{49/2}.
\end{eqnarray*}
Recall that $h^4-4c^2$ is strictly positive. By setting $L_1=0$, we seek the solutions of
$-h c^{7/2} \left(h^4+4 c^2\right)^2 =0.$
Thus, the only possible conditions are $h=0$ or $c=0$. However, neither condition is possible, since $c>0$ and $h>0$. Hence, $L_1\neq0$ and so $E_4^-$ is not a center on the center manifold. In summary, under the conditions $a=-c$, $b>3c$ and $c>0$, $E_4^-$ is a focus. So, $E_4^+$ is also a focus on the center manifold.

Now, let us consider points $E_5^\pm$. The computations here are very similar to those made at points $E_4^\pm$. Likewise, we will consider point $E_5^-$ for the computations due to the symmetry of the points.

The characteristic polynomial of Jacobian matrix $J(E_5^-)$  is $P(\lambda)=\lambda^3+\alpha\lambda^2+\beta \lambda+\gamma$, where
\begin{eqnarray*}
	\alpha&=&a-c;\\
	\beta&=&-\dfrac{a^2+a b +a c + b c+a\sqrt{(a+b)^2+4a c}-c\sqrt{(a+b)^2+4a c}}{2a c};\\
	\gamma&=&2\sqrt{(a+b)^2+4a c}.
\end{eqnarray*}
The conditions for $E_5^-$ to be a Hopf point are $a=-c$, $b<0$, $c<0$ and $b<3c$. The eigenvalues are $\lambda_{1,2}=\pm\dfrac{\sqrt[4]{(b-3c)(b+c)}}{\sqrt{-c}}i$ and $\lambda_3=2c$. Now, translating point $E_5^-$ to the origin, inserting a new real parameter $h>0$ given by the equation
\begin{equation*}
	\sqrt{(b-3 c) (b+c)}-b+c-h^2=0,
\end{equation*}
that is,
\begin{equation*}
	b=\frac{-4 c^2+2 c h^2-h^4}{2 h^2},
\end{equation*}
and making the following change of coordinates
\begin{equation*}
	\begin{aligned}
		(x,y,z)\mapsto&\Bigg(\dfrac{2 c \left(c h^2-1\right) \sqrt{h^4-4 c^2}}{8 c^3 h^2+4 c^2-h^4}u-\dfrac{ \sqrt{-c} h\left(4 c^2+h^4\right)}{\sqrt{2} \left(8 c^3 h^2+4 c^2-h^4\right)}v-\dfrac{h^2}{2 c}w,\\
		&\quad\dfrac{ \left(4 c^3-h^2\right) \sqrt{h^4-4 c^2}}{8 c^3 h^2+4 c^2-h^4}u-\dfrac{\sqrt{2} (-c)^{3/2} \left(4 c^2+h^4\right)}{8 c^3 h^3+4 c^2 h-h^5}v+w, v\Bigg)
	\end{aligned}
\end{equation*}
with the time rescaling $\tau=\dfrac{\sqrt{h^4-4 c^2} }{\sqrt{2} \sqrt{-c} h}t$, we can write the system \eqref{sist3Dchaotic} as follows

\begin{equation}\label{sistE5mnorm}
	\begin{aligned}
		\dot{u}={}&-v+\dfrac{h \left(4 c^2+h^4\right)}{\sqrt{2} \sqrt{-c} \left(4 c^2-h^4\right)}v w-\dfrac{c \left(4 c^2+h^4\right)^2}{\left(4 c^2-h^4\right) \left(8 c^3 h^2+4 c^2-h^4\right)}v^2+\dfrac{2 \sqrt{2} c \sqrt{-c} h \left(-4 c^3+c h^4+2 h^2\right)}{\sqrt{h^4-4 c^2} \left(8 c^3 h^2+4 c^2-h^4\right)}u v\\
		\dot{v}={}& u-\dfrac{h^3 }{\sqrt{2} \sqrt{-c} \sqrt{h^4-4 c^2}}w^2-\dfrac{c \left(4 c^2+h^4\right)^2}{\left(8 c^3 h^2+4 c^2-h^4\right)^2}u v+\dfrac{h \left(4 c^2+h^4\right)}{\sqrt{2} \sqrt{-c} \left(8 c^3 h^2+4 c^2-h^4\right)}u w\\
		& -\dfrac{2 c h^2\left(4 c^2+h^4\right)}{\sqrt{h^4-4 c^2} \left(8 c^3 h^2+4 c^2-h^4\right)}v w-\dfrac{2 \sqrt{2} c\sqrt{-c} h \left(4c^3-h^2\right) \left(1+c h^2\right) \sqrt{h^4-4 c^2}}{\left(8 c^3 h^2+4 c^2-h^4\right)^2}u^2\\
		& +\dfrac{\sqrt{2} c^{3} h \left(4 c^2+h^4\right)^2}{\sqrt{-c}\sqrt{h^4-4 c^2} \left(8 c^3 h^2+4 c^2-h^4\right)^2}v^2\\
		\dot{w}={}&\dfrac{2 \sqrt{2} (-c)^{3/2} h }{\sqrt{h^4-4 c^2}}w+\dfrac{4 c^3 \left(2 c-h^2\right) \left(2 c+h^2\right) \left(4 c^3-h^2\right) \left(c h^2+1\right) \left(4 c^2+h^4\right)}{\sqrt{h^4-4 c^2} \left(8 c^3 h^2+4 c^2-h^4\right)^3}u^2 \\
		&+\dfrac{2 c^3 \left(4 c^2+h^4\right) \left(64 c^6 h^2+48 c^5-16 c^4 h^6-40 c^3 h^4-16 c^2 h^2+3 c h^8+4 h^6\right)}{\sqrt{h^4-4 c^2} \left(8 c^3 h^2+4 c^2-h^4\right)^3}v^2\\
		&+\dfrac{c h^2 \left(4 c^2+h^4\right)}{\sqrt{h^4-4 c^2} \left(8 c^3 h^2+4 c^2-h^4\right)} w^2+\dfrac{4 \sqrt{2} c\sqrt{-c} h \left(2 c-h^2\right) \left(2 c+h^2\right) \left(4 c^3-h^2\right) \left(c h^2+1\right)}{\sqrt{h^4-4 c^2} \left(8 c^3 h^2+4 c^2-h^4\right)^2} v w\\
		&-\dfrac{\sqrt{2} (-c)^{3/2} \left(4 c^2+h^4\right) \left(64 c^7 h^4+32 c^6 h^2-16 c^5-8 c^4 h^6+8 c^3 h^4+8 c^2 h^2-c h^8-2 h^6\right)}{h \left(8 c^3 h^2+4 c^2-h^4\right)^3} u v\\
		&-\frac{64 c^7 + 128 c^8 h^2 - 16 c^4 \sqrt{-c} h^4 + 64 c^6 h^6 - 4 c (-c)^{3/2} h^8 + 8 c^4 h^{10} - c h^{12}}{\left(8 c^3 h^2 + 4 c^2 - h^4\right)^3}u w
	\end{aligned}
\end{equation}

Using the algorithm of Section \ref{secCenter} for obtaining focus quantities, we have
\begin{equation*}
	L_1=h (-c)^{7/2} \sqrt{h^4-4 c^2} \left(h^4+4 c^2\right)^2.
\end{equation*}
Again, as $h^4-4c^2>0$, setting $L_1=0$, the solutions come from $h(-c)^{7/2} \left(h^4+4 c^2\right)^2=0$. Hence, the solutions are $c=0$ or $h=0$. But, none of these solutions are possible, because $c<0$ and $h>0$. Therefore, $E_5^-$, and consequently $E_5^+$, are foci on the center manifold. \qed

\section{Numerical simulations}\label{secSimul}
In this section, we will show some simulations that illustrate the results of Theorems \ref{Teo1} and \ref{Teo2} and Corollary \ref{coro}.

\subsection{Simulations of Theorem \ref{Teo1}}

Figure~\ref{fig1} shows the orbits near the origin of the three-dimensional differential system \eqref{sistE1cond6} for \( d = 1 \) (i.e., the origin is a center on the center manifold), under different initial conditions. Although the explicit expression of the center manifold is unknown, Figure~\ref{fig1} provides insight into its shape by illustrating how the orbits of the system evolve along it. Figure~\ref{fig2} shows the time evolution of each variable for the initial condition \((0.5, -0.75, 0.1)\), which is close to the origin. This allows for a clearer observation of the periodic behavior of the system.

\begin{figure}[h!]
	\centering
	
	\begin{subfigure}[b]{0.49\textwidth}
		\includegraphics[width=\linewidth]{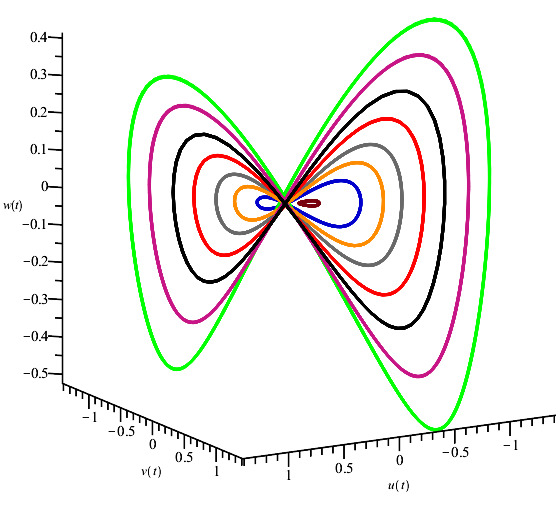}
		%\caption{Periodic behavior of $x(t)$}
	\end{subfigure}
	\hfill
	\begin{subfigure}[b]{0.49\textwidth}
		\includegraphics[width=\linewidth]{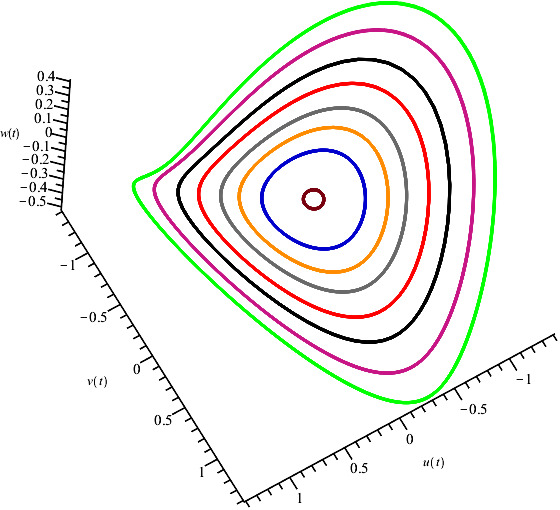}
		%\caption{Periodic behavior of $y(t)$}
	\end{subfigure}
	
%	\vspace{\baselineskip}
	
	\begin{subfigure}[b]{0.49\textwidth}
		\includegraphics[width=\linewidth]{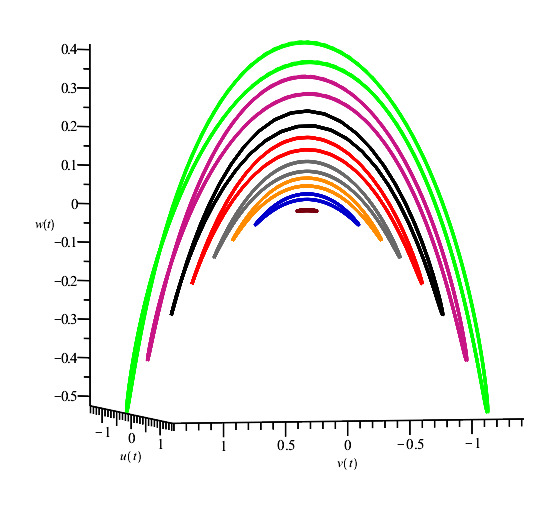}
		%\caption{Periodic behavior of $z(t)$}
	\end{subfigure}
	\hfill
	\begin{subfigure}[b]{0.49\textwidth}
		\includegraphics[width=\linewidth]{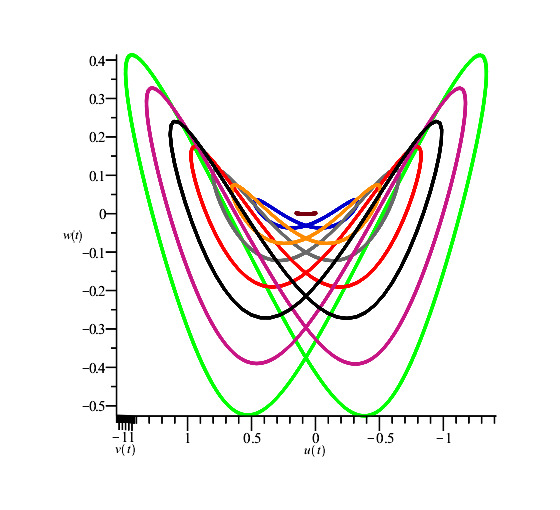}
		%\caption{Periodic behavior of $x(t), y(t)$ and $z(t)$}
	\end{subfigure}
	
	\caption{Orbits of system \eqref{sistE1cond6} in $uvw$-space with initial conditions $(0.08,0.002,0.03), (0.4,0.07,0.13)$,$(-0.5,0.3,0.25)$,$(0.2,0.7,0.85)$, $(0.5,0.75,0.5)$,$(0.8,0.7,-0.5)$,$(-1,-0.75,0.6)$ and $(-1,1,1)$, for $d=1$. }
	\label{fig1}
\end{figure}

\begin{figure}[h!]
	\centering
	
	\begin{subfigure}[b]{0.45\textwidth}
		\includegraphics[width=\linewidth]{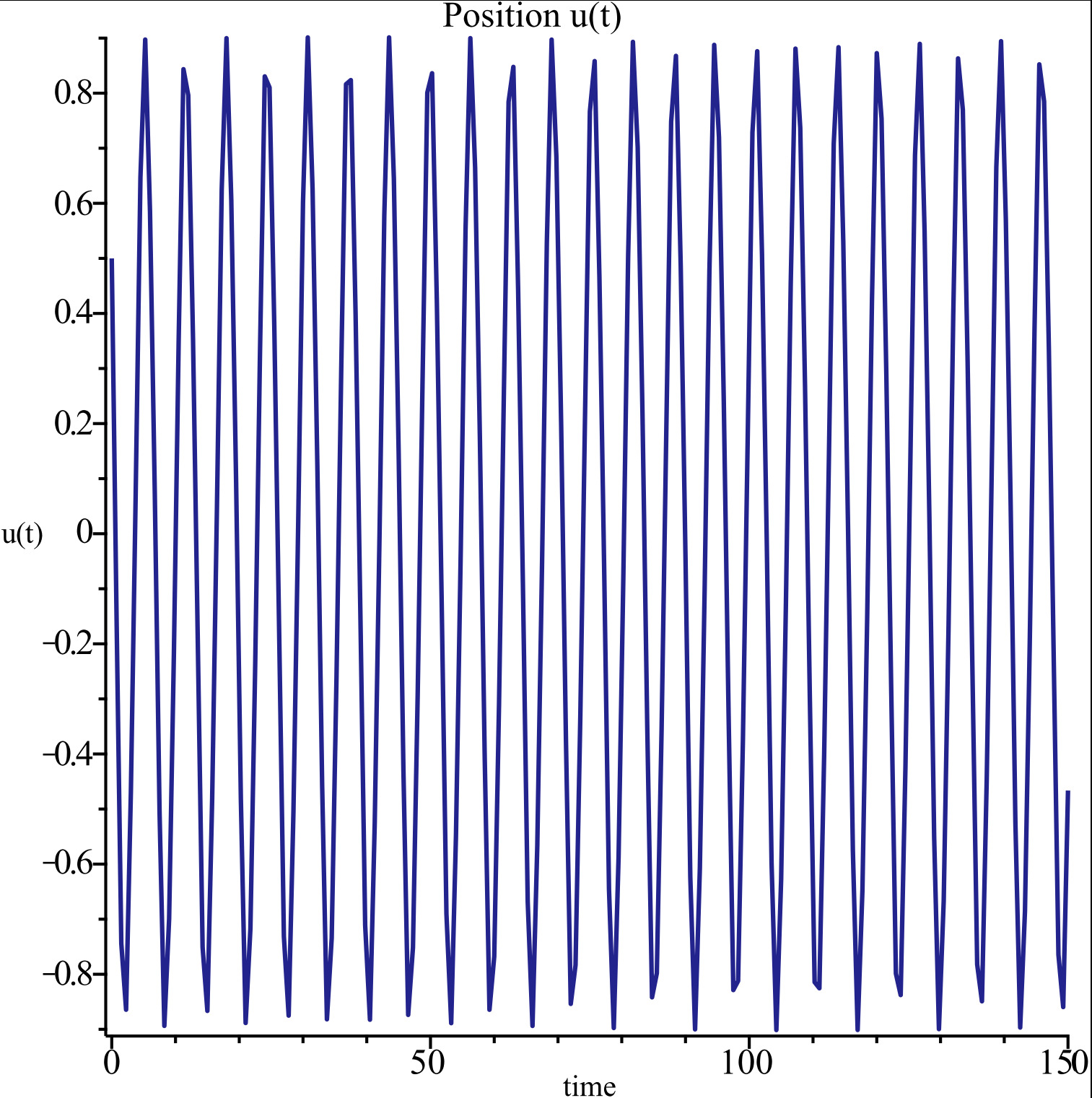}
		\caption{Periodic behavior of $u(t)$}
		\label{fig:sub_a}
	\end{subfigure}
	\hfill
	\begin{subfigure}[b]{0.45\textwidth}
		\includegraphics[width=\linewidth]{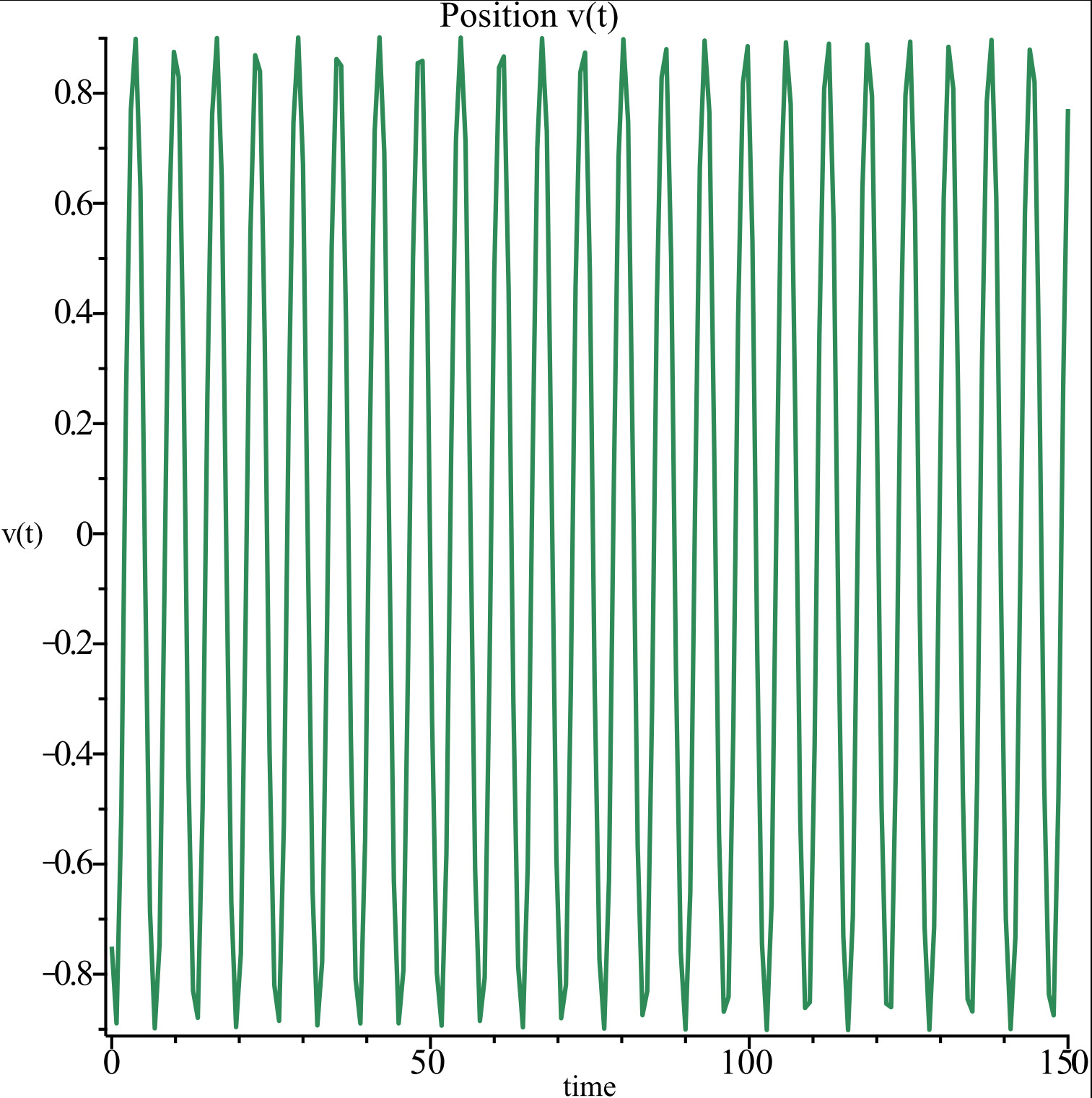}
		\caption{Periodic behavior of $v(t)$}
		\label{fig:sub_b}
	\end{subfigure}
	
%	\vspace{\baselineskip}
	
	\begin{subfigure}[b]{0.45\textwidth}
		\includegraphics[width=\linewidth]{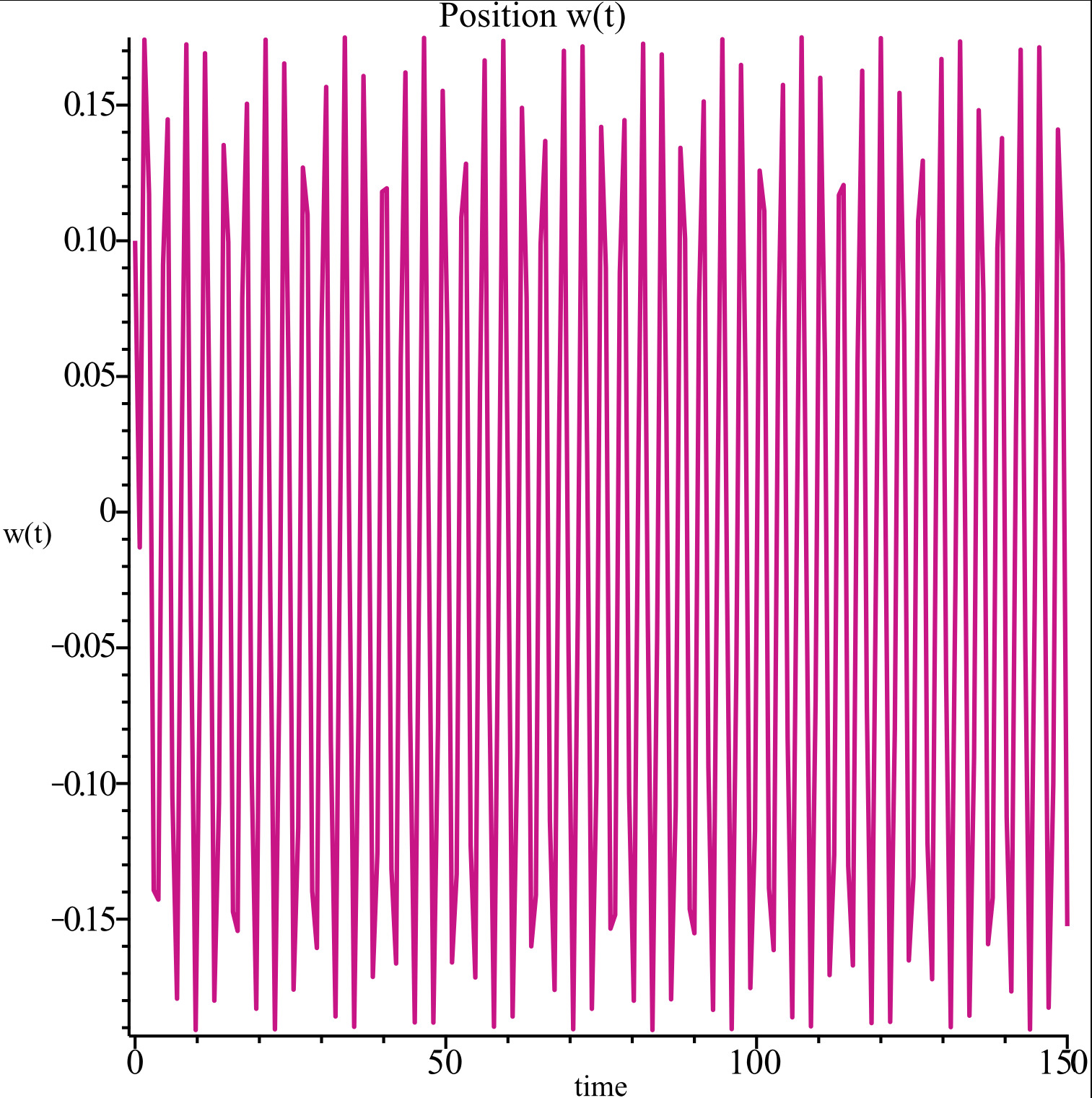}
		\caption{Periodic behavior of $w(t)$}
		\label{fig:sub_c}
	\end{subfigure}
	\hfill
	\begin{subfigure}[b]{0.45\textwidth}
		\includegraphics[width=\linewidth]{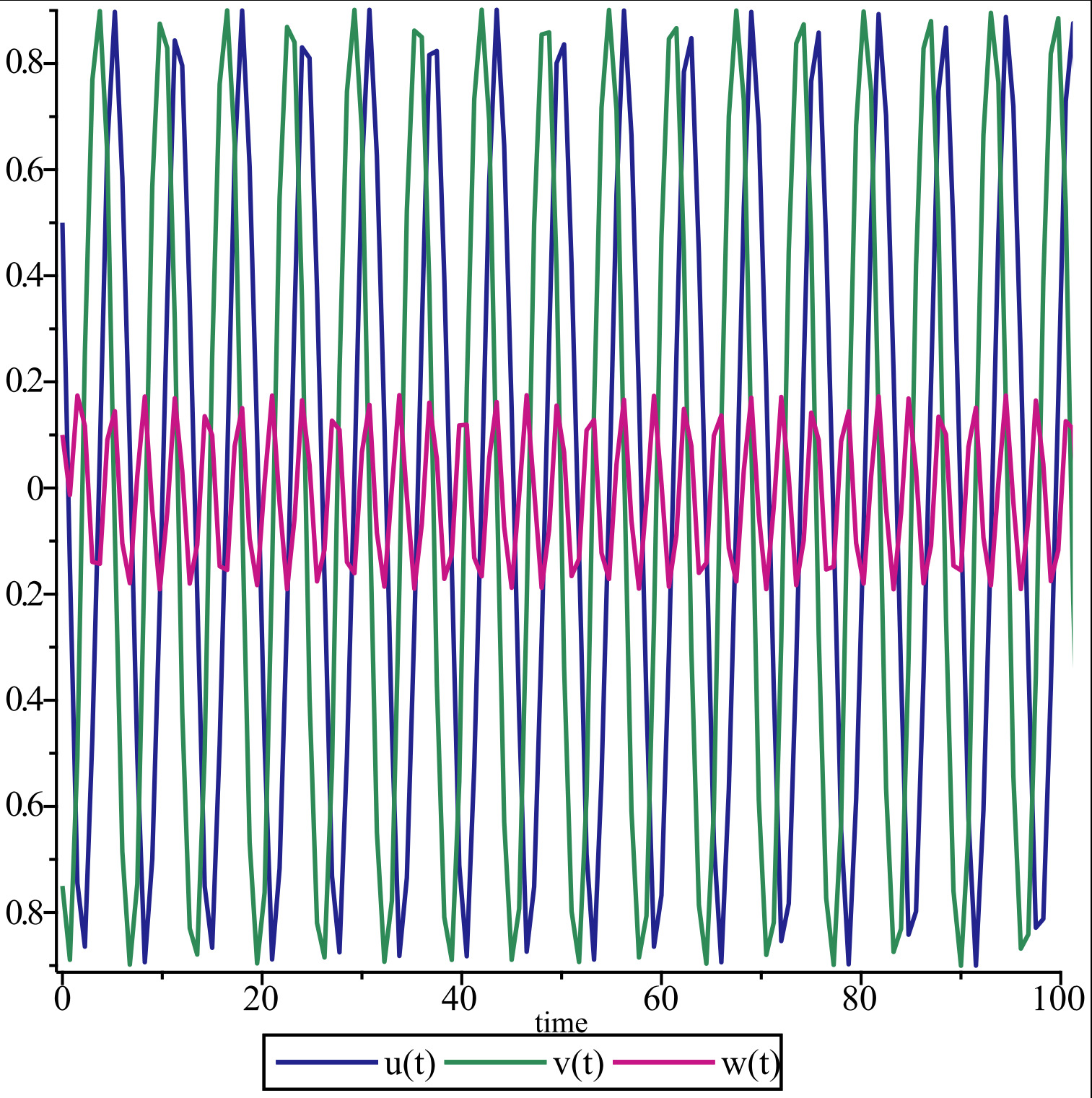}
		\caption{Periodic behavior of $u(t), v(t)$ and $w(t)$}
		\label{fig:sub_d}
	\end{subfigure}
	
	\caption{Behavior of each variable over time for initial condition $(0.5,-0.75,0.1)$}
	\label{fig2}
\end{figure}

\subsection{Simulations of Corollary \ref{coro}}
Figures \ref{figE2} and \ref{figE3} illustrate the dynamics near the focus points $E_2^\pm$ and $E_3^\pm$ under the conditions presented in Corollary \ref{coro}. The simulations are presented for both forward and backward time, considering several initial conditions.

\begin{figure}[h!]
\centering

\begin{subfigure}[b]{0.49\textwidth}
	\includegraphics[width=\linewidth]{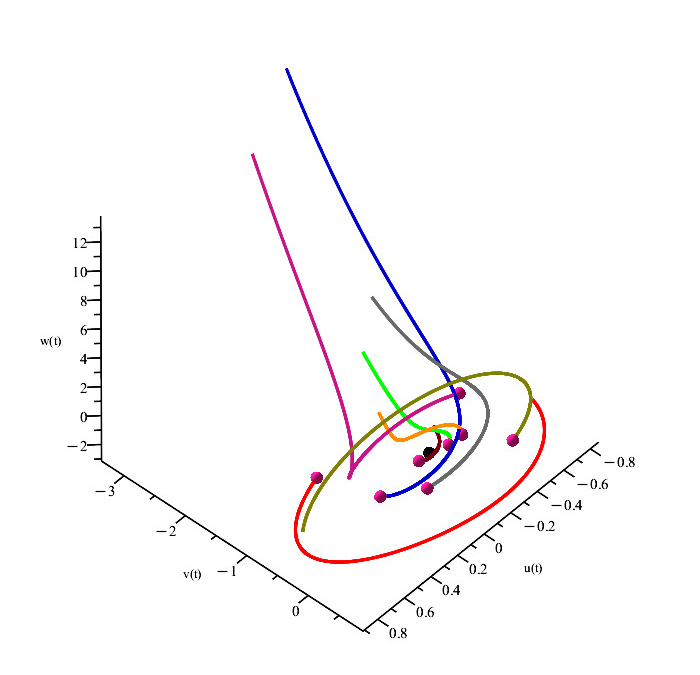}
	\caption{Forward time.}
\end{subfigure}
\hfill
\begin{subfigure}[b]{0.49\textwidth}
	\includegraphics[width=\linewidth]{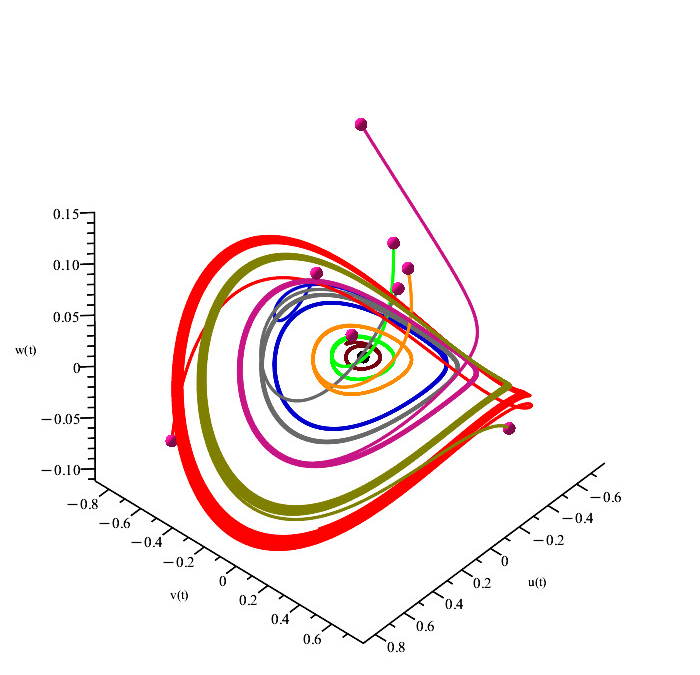}
	\caption{Backward time.}
\end{subfigure}

\caption{Orbits of system  \eqref{sistE2norm} in $uvw$-space with initial conditions $(0.4, 0.07, 0.13)$, $(0.08, 0.002, 0.03)$, $(-0.1, 0.1, 0.11)$, $(0.2, 0.4, 0.125)$, $(0.5, -0.375, -0.1)$, $(-0.2, 0.1, 0.075)$, $(-0.6, -0.375, 0.15)$ and $(-0.35, 0.6, -0.05)$.}
\label{figE2}
\end{figure}

\begin{figure}[h!]
\centering

\begin{subfigure}[b]{0.49\textwidth}
	\includegraphics[width=\linewidth]{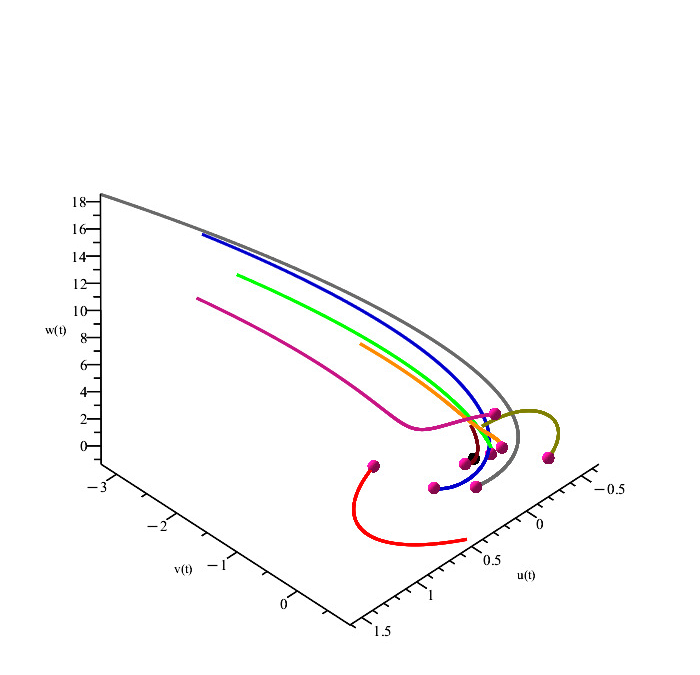}
	\caption{Forward time.}
\end{subfigure}
\hfill
\begin{subfigure}[b]{0.49\textwidth}
	\includegraphics[width=\linewidth]{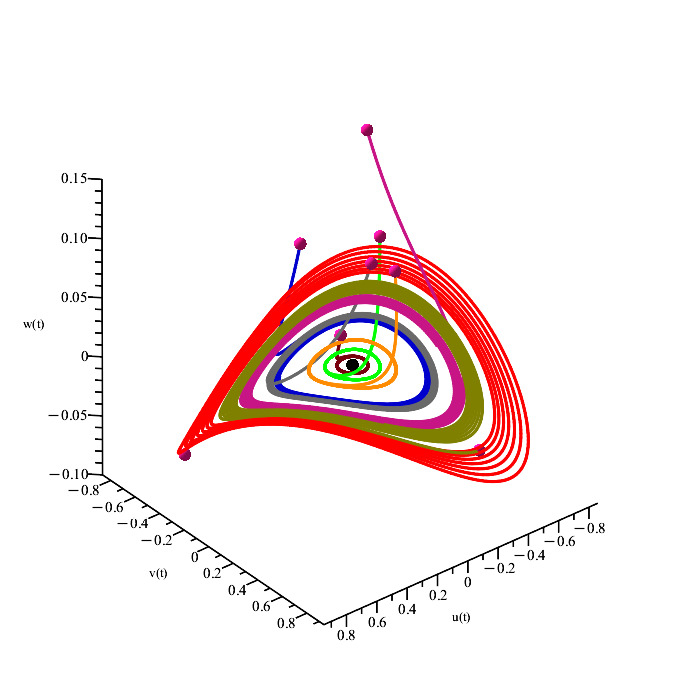}
	\caption{Backward time.}
\end{subfigure}

\caption{Orbits of system  \eqref{sistE3norm} in $uvw$-space with initial conditions $(0.4, 0.07, 0.13)$, $(0.08, 0.002, 0.03)$, $(-0.1, 0.1, 0.11)$, $(0.2, 0.4, 0.125)$, $(0.5, -0.375, -0.1)$, $(-0.2, 0.1, 0.075)$, $(-0.6, -0.375, 0.15)$ and $(-0.35, 0.6, -0.05)$.}
\label{figE3}
\end{figure}

\subsection{Simulations of Theorem \ref{Teo2}}
Figures~\ref{fig3} and~\ref{fig4} show simulations of the behavior of the orbits of system~\eqref{sistE4mnorm} near the origin, which in this case is a focus on the center manifold. The simulations are presented for both forward and backward time, considering different values of the system parameters and several initial conditions. Figures~\ref{fig5} and~\ref{fig6} present similar simulations for the orbits of system~\eqref{sistE5mnorm}.
\begin{figure}[h!]
	\centering
	
	\begin{subfigure}[b]{0.49\textwidth}
		\includegraphics[width=\linewidth]{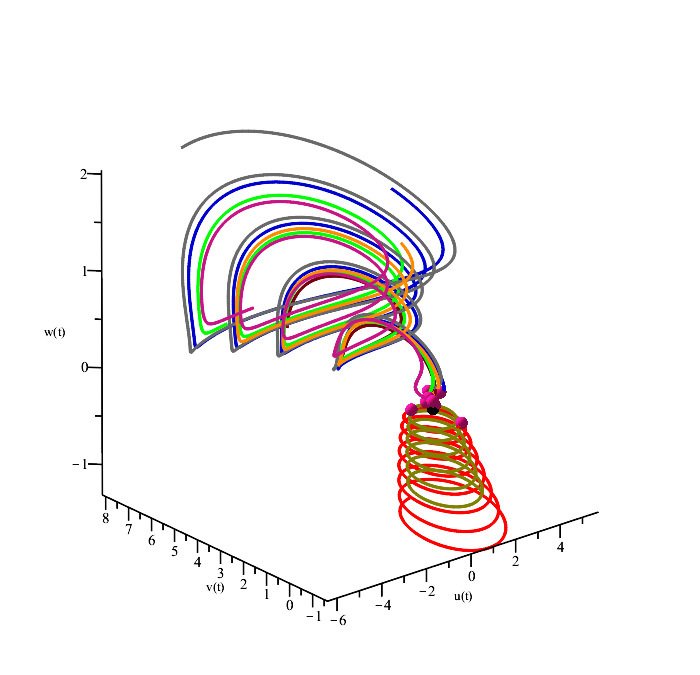}
		\caption{Forward time for $c=\frac{1}{4}$ and $h=2$.}
	\end{subfigure}
	\hfill
	\begin{subfigure}[b]{0.49\textwidth}
		\includegraphics[width=\linewidth]{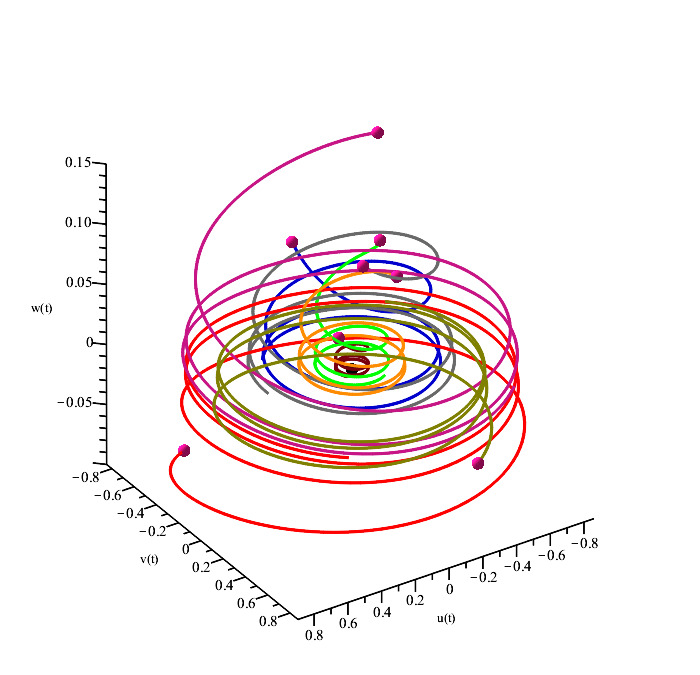}
		\caption{Backward time for $c=\frac{1}{4}$ and $h=2$.}
	\end{subfigure}
	
%	\vspace{\baselineskip}
	
	\begin{subfigure}[b]{0.49\textwidth}
		\includegraphics[width=\linewidth]{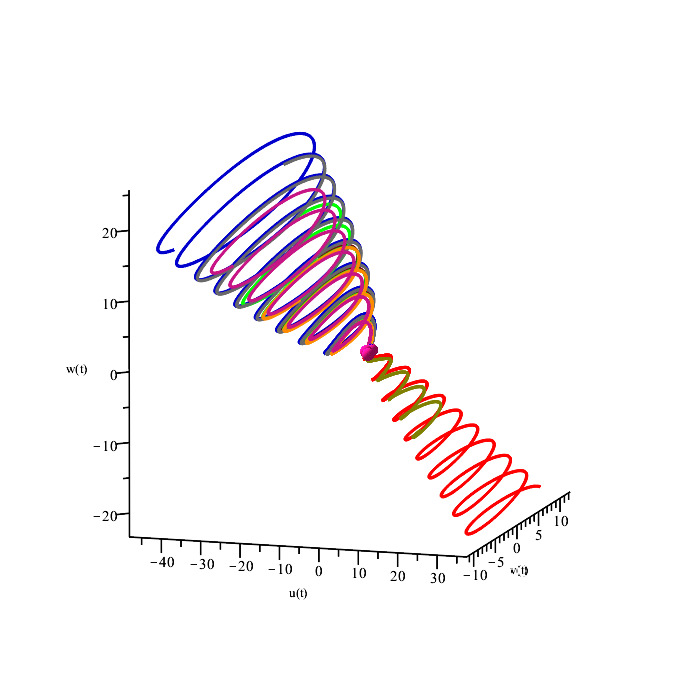}
		\caption{Forward time for $c=1$ and $h=2$.}
	\end{subfigure}
	\hfill
	\begin{subfigure}[b]{0.49\textwidth}
		\includegraphics[width=\linewidth]{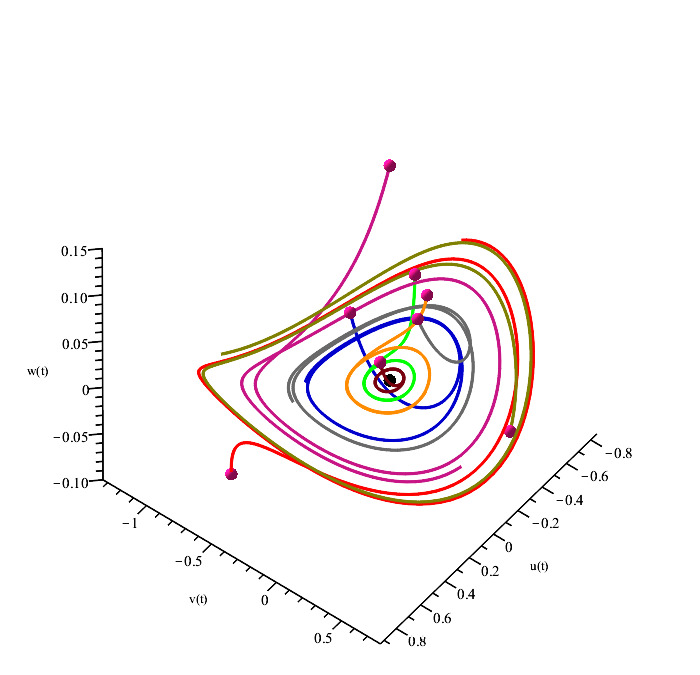}
		\caption{Backward time for $c=1$ and $h=2$.}
	\end{subfigure}
	
	\caption{Orbits of system  \eqref{sistE4mnorm} in $uvw$-space with initial conditions $(0.4, 0.07, 0.13)$, $(0.08, 0.002, 0.03)$, $(-0.1, 0.1, 0.11)$, $(0.2, 0.4, 0.125)$, $(0.5, -0.375, -0.1)$, $(-0.2, 0.1, 0.075)$, $(-0.6, -0.375, 0.15)$ and $(-0.35, 0.6, -0.05)$. For $c=\frac{1}{4}$ and $h=2$, the first focus quantity is $L_1=-\frac{3}{2032}\sqrt{\frac{7}{2}}$. For $c=1$ and $h=2$, we obtain $L_1=-\frac{5}{121}\sqrt{\frac{3}{2}}$.}
	\label{fig3}
\end{figure}

\begin{figure}[h!]
	\centering
	
	\begin{subfigure}[b]{0.49\textwidth}
		\includegraphics[width=\linewidth]{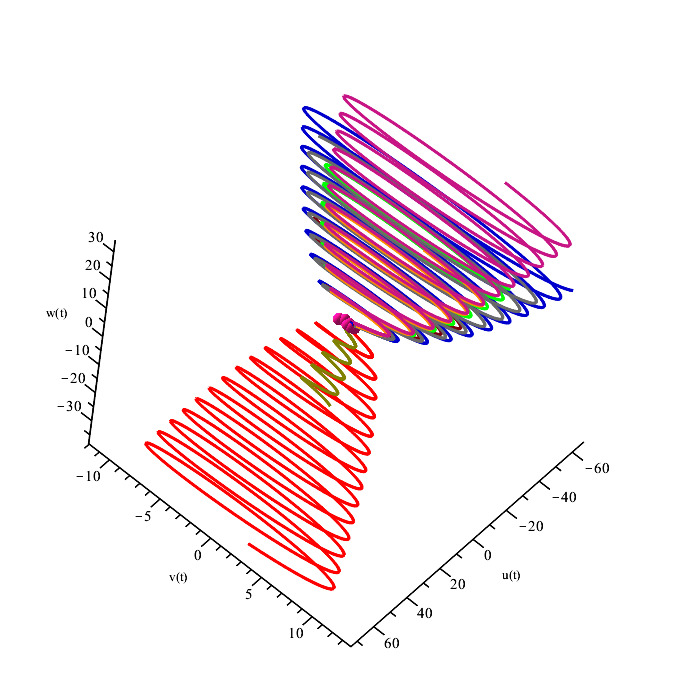}
		\caption{Forward time for $c=3$ and $h=5$.}
	\end{subfigure}
	\hfill
	\begin{subfigure}[b]{0.49\textwidth}
		\includegraphics[width=\linewidth]{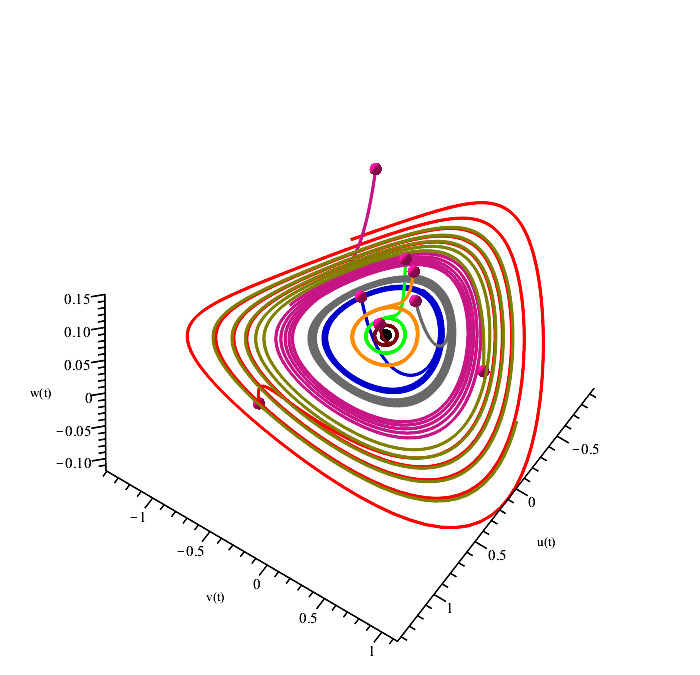}
		\caption{Backward time for $c=3$ and $h=5$.}
	\end{subfigure}
	
%	\vspace{\baselineskip}
	
	\begin{subfigure}[b]{0.49\textwidth}
		\includegraphics[width=\linewidth]{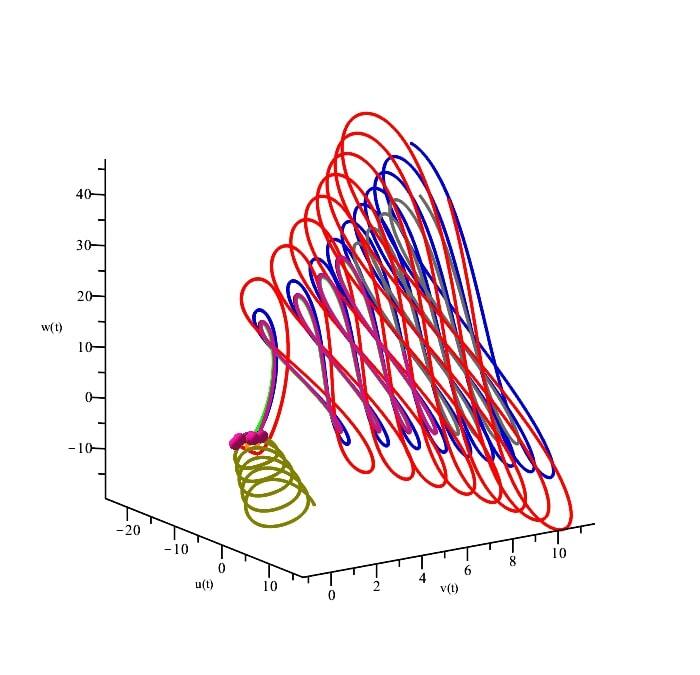}
		\caption{Forward time for $c=\frac{1}{9}$ and $h=\frac{1}{2}$.}
	\end{subfigure}
	\hfill
	\begin{subfigure}[b]{0.49\textwidth}
		\includegraphics[width=\linewidth]{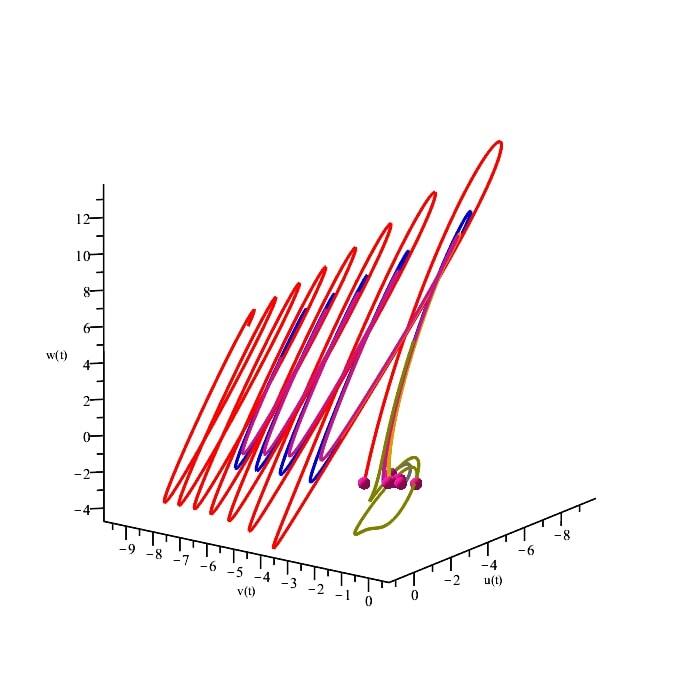}
		\caption{Backward time for $c=\frac{1}{9}$ and $h=\frac{1}{2}$.}
	\end{subfigure}
	
	\caption{Orbits of system  \eqref{sistE4mnorm} in $uvw$-space with initial conditions $(0.4, 0.07, 0.13)$, $(0.08, 0.002, 0.03)$, $(-0.1, 0.1, 0.11)$, $(0.2, 0.4, 0.125)$, $(0.5, -0.375, -0.1)$, $(-0.2, 0.1, 0.075)$, $(-0.6, -0.375, 0.15)$ and $(-0.35, 0.6, -0.05)$. For $c=3$ and $h=5$, the first focus quantity is $L_1=-\frac{58984335 \sqrt{\frac{1767}{2}}}{69548286619}$. For $c=1$ and $h=2$, we obtain $L_1=-\frac{2523 \sqrt{34}}{220409}$. 
	} 
	\label{fig4}
\end{figure}

\begin{figure}[h!]
	\centering
	
	\begin{subfigure}[b]{0.49\textwidth}
		\includegraphics[width=\linewidth]{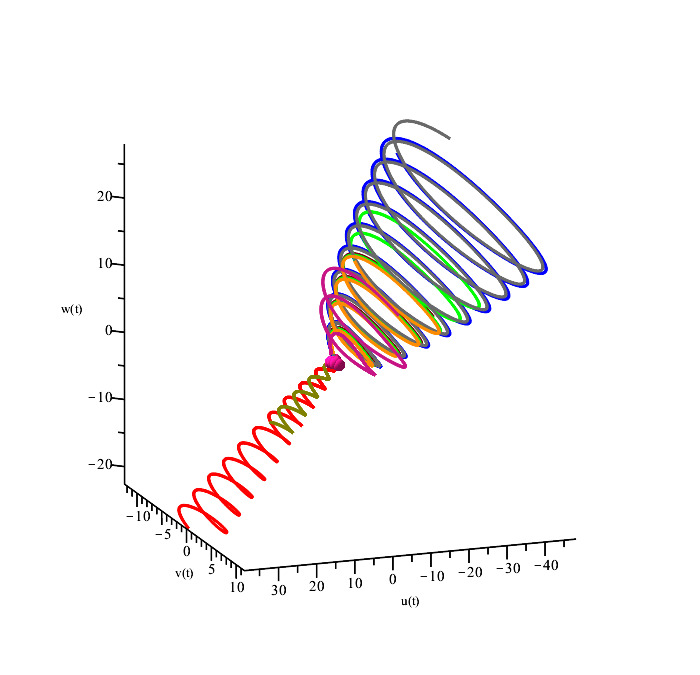}
		\caption{Forward time for $c=-1$ and $h=2$.}
	\end{subfigure}
	\hfill
	\begin{subfigure}[b]{0.49\textwidth}
		\includegraphics[width=\linewidth]{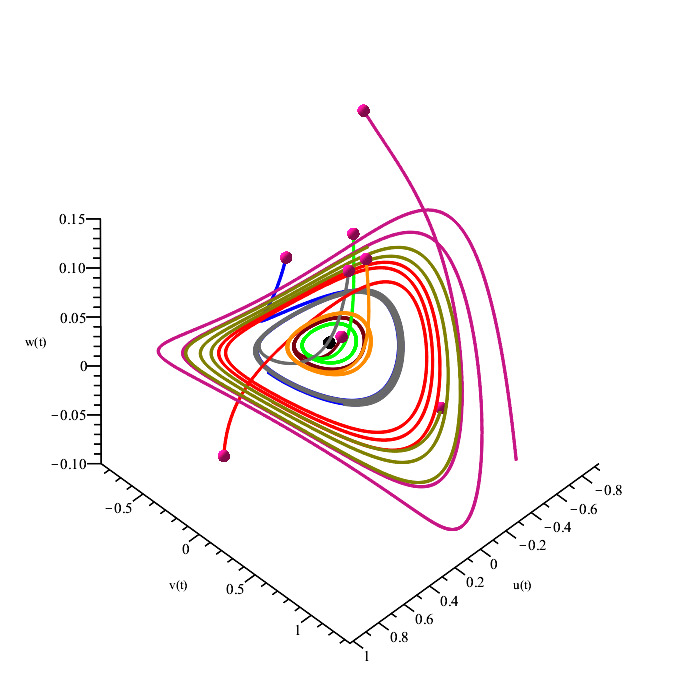}
		\caption{Backward time for $c=-1$ and $h=2$.}
	\end{subfigure}
	
	%	\vspace{\baselineskip}
	
	\begin{subfigure}[b]{0.49\textwidth}
		\includegraphics[width=\linewidth]{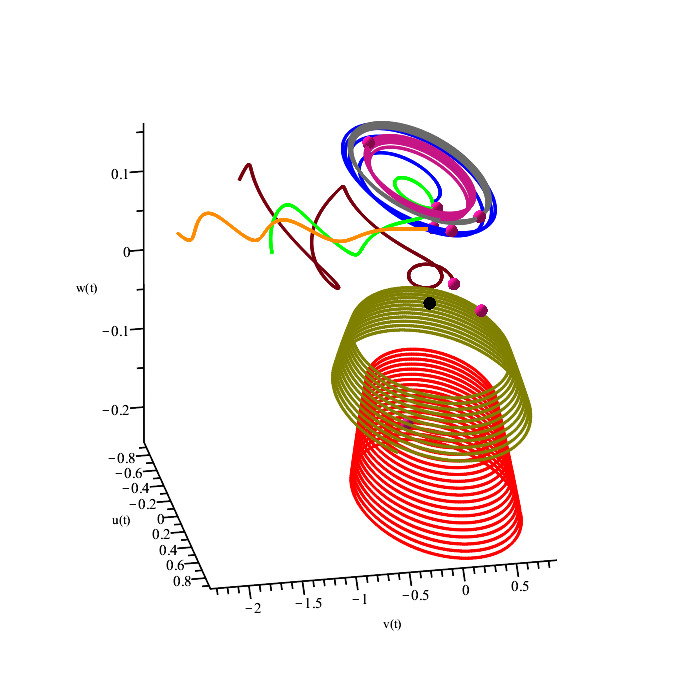}
		\caption{Forward time for $c=-\frac{1}{4}$ and $h=10$.}
	\end{subfigure}
	\hfill
	\begin{subfigure}[b]{0.49\textwidth}
		\includegraphics[width=\linewidth]{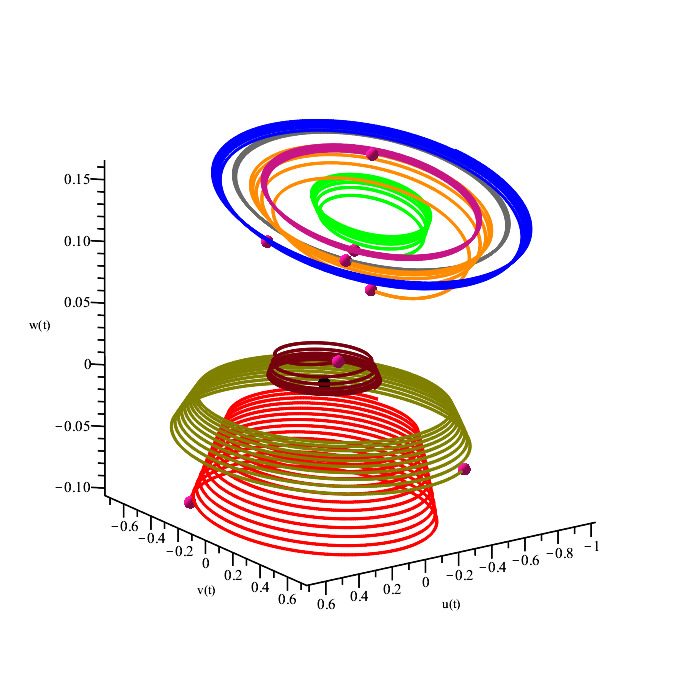}
		\caption{Backward time for $c=-\frac{1}{4}$ and $h=10$.}
	\end{subfigure}
	
	\caption{Orbits of system  \eqref{sistE5mnorm} in $uvw$-space with initial conditions $(0.4, 0.07, 0.13)$, $(0.08, 0.2, 0.03)$, $(-0.1, 0.1, 0.11)$, $(0.2, 0.4, 0.125)$, $(0.5, -0.375, -0.1)$, $(-0.2, 0.1, 0.075)$, $(-0.6, -0.375, 0.15)$ and $(-0.35, 0.6, -0.05)$. For $c=-1$ and $h=2$, the first focus quantity is $L_1=-\frac{5 \sqrt{\frac{3}{2}}}{121}$. For $c=-\frac{1}{4}$ and $h=10$, we obtain $L_1=-\frac{8000400005 \sqrt{\frac{39999}{2}}}{2053610916723568}$.
	}
	\label{fig5}
\end{figure}
\begin{figure}[h!]
	\centering
	
	\begin{subfigure}[b]{0.49\textwidth}
		\includegraphics[width=\linewidth]{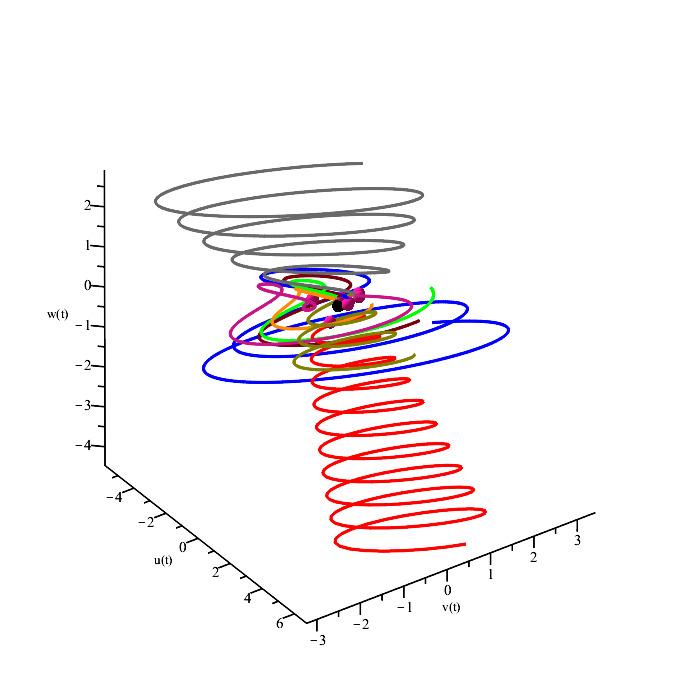}
		\caption{Forward time for $c=-2$ and $h=10$.}
	\end{subfigure}
	\hfill
	\begin{subfigure}[b]{0.49\textwidth}
		\includegraphics[width=\linewidth]{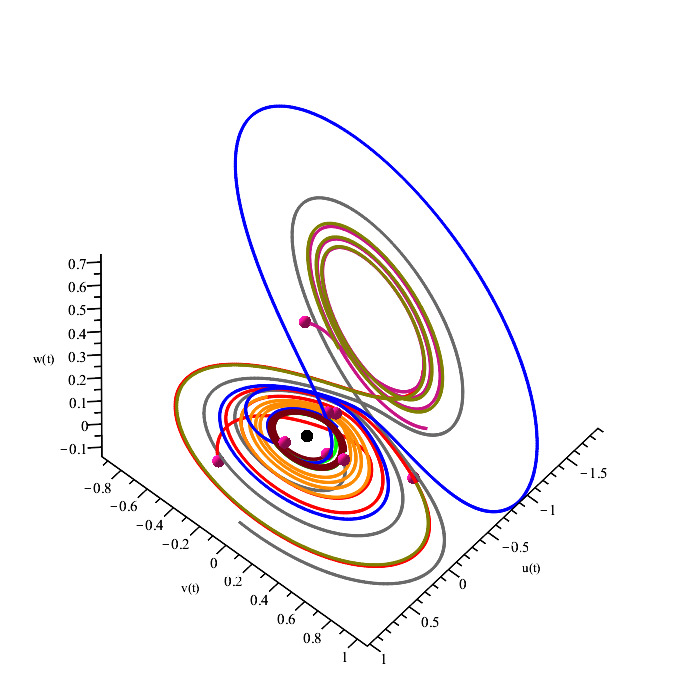}
		\caption{Backward time for $c=-2$ and $h=10$.}
	\end{subfigure}
	
	%	\vspace{\baselineskip}
	
	\begin{subfigure}[b]{0.49\textwidth}
		\includegraphics[width=\linewidth]{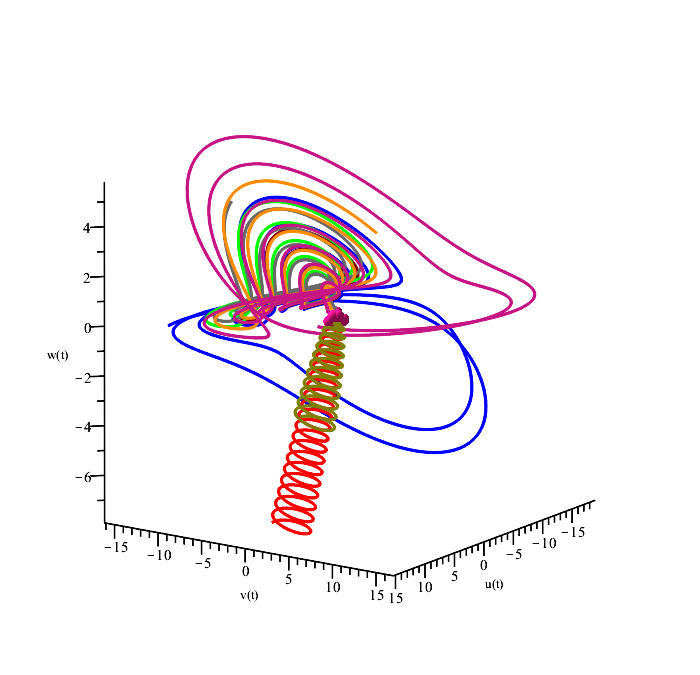}
		\caption{Forward time for $c=-\frac{1}{3}$ and $h=2$.}
	\end{subfigure}
	\hfill
	\begin{subfigure}[b]{0.49\textwidth}
		\includegraphics[width=\linewidth]{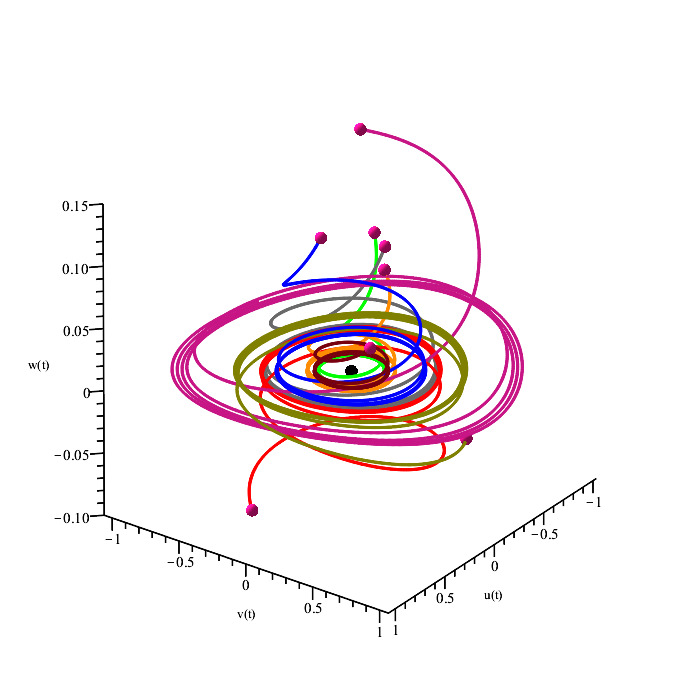}
		\caption{Backward time for $c=-\frac{1}{3}$ and $h=2$.}
	\end{subfigure}
	
	\caption{Orbits of system  \eqref{sistE5mnorm} in $uvw$-space with initial conditions $(0.4, 0.07, 0.13)$, $(0.08, 0.2, 0.03)$, $(-0.1, 0.1, 0.11)$, $(0.2, 0.4, 0.125)$, $(0.5, -0.375, -0.1)$, $(-0.2, 0.1, 0.075)$, $(-0.6, -0.375, 0.15)$ and $(-0.35, 0.6, -0.05)$. For $c=-2$ and $h=10$, the first focus quantity is $L_1=-\frac{489845 \sqrt{39}}{11862016}$. For $c=-\frac{1}{3}$ and $h=2$, we obtain $L_1=-\frac{1369 \sqrt{\frac{105}{2}}}{1366283}$.
	}
	\label{fig6}
\end{figure}

\section{Cyclicity}\label{secCyc}
The study of the cyclicity of a Hopf point in three-dimensional differential systems is closely linked to computing the Lyapunov coefficients (or, equivalently, the focus quantities), since these are the coefficients of the power series expansion of the reduced displacement function, whose isolated zeros correspond to limit cycles. The author of \cite{Chris05}, in the planar case, explores this relationship. This work analyzes how only the linear parts of the Lyapunov coefficients can be used to estimate a lower bound for the number of limit cycles that bifurcate from a Hopf point of the system. There are also results that use the higher-order parts of Lyapunov coefficients to detect additional limit cycles, as well as efficient methods that enable the computation of these parts at a lower cost, from a computational point of view (see \cite{GouveiaTorre21}).

Let us remember that, in the parameter space, the center conditions constitute the so-called \textit{Bautin variety}. In fact, a point $s$ on the Bautin variety corresponds to a point satisfying $L_k(s) = 0$, for all $k \geq 1$.

The following results were demonstrated in \cite{GouveiaTorre21} (see also \cite{teselucas,Chris05,tesegouveia,Han99,Pessoa24}).

%A fairly common method for studying the cyclicity of a Hopf point in a three-dimensional differential system is to compute the focus quantities (or the Lyapunov coefficients, equivalently) and study the signs, as well as the zeros, of the reduced displacement function from small perturbations of the parameters, since isolated zeros correspond to limit cycles. This idea was initially formulated for planar systems. In \cite{Chris05}, it was proved that, by analyzing only the linear parts of the focus quantities, it is possible to obtain an estimate for this number of limit cycles. Furthermore higher orders of the focus quantities are used to obtain more limit cycles, and the parallelization method, extended in \cite{GouveiaTorre21} to higher orders, is an extremely useful tool for this. The following results were demonstrated in \cite{Chris05}. But first, we call the set of points in the parameter space that makes the system have a center at the origin on the center manifold the Bautin variety.  In fact, a point on the Bautin manifold is what we call the center condition and we can also view this as a point satisfying $L_k(s)=0$, for all $k\geq1$.
\begin{theorem}\label{TeoCyc1}
	Suppose that $s$ is a point on the Bautin variety and that the first $k$ focus quantities, $L_1, \dots, L_k$, have independent linear parts (with respect to the expansion of $L_i$ about $s$). Then $s$ lies on a component of the Bautin variety of codimension at least $k$ and there are bifurcations which produce $k$ limit cycles locally from the center corresponding to the parameter value $s$. If, furthermore, we know that $s$ lies on a component of the center variety of codimension $k$, then $s$ is a smooth point of the variety, and the cyclicity of the center for the parameter value $s$ is exactly $k$. In the latter case, $k$ is also the cyclicity of a generic point on this component of the Bautin variety.
\end{theorem} 

\begin{theorem}\label{TeoCyc2}
	Suppose that $s$ is a point in the Bautin variety such that the $k$ first focus quantities have independent linear parts (with respect to the expansion of $L_i$ about $s$) and the next $l$ ones have their linear parts as linear combinations of the linear parts of the $k$ first ones. After a change of variables in the parameter space, if necessary, we can write $L_i = u_i$, for $i = 1, \dots, k$ and, assuming $L_0 = L_1 = \cdots = L_k = 0$, the next focus quantities are $L_i = h_i(u) + O(|u|^{m+1})$, for $i = k+1, \dots, k+l$, where $h_i$ are homogeneous polynomials of degree $m \geq 2$ and $u = (u_{k+1}, \dots, u_{k+l})$. If there exists a line $\eta$ in the parameter space such that for $k < i < k+l$, $h_i(\eta) = 0$, the hypersurfaces $h_i = 0$ intersect transversally along $\eta$, and $h_{k+l}(\eta) \neq 0$, then there are perturbations of the center which produce $k+l$ limit cycles.
\end{theorem}

Theorems \ref{TeoCyc1} and \ref{TeoCyc2} found in \cite{GouveiaTorre21} were originally stated for the Lyapunov coefficients, but by virtue of Theorem 7 and Proposition 34 of \cite{Gar18}, we adapted the above statements for focus quantities (see also \cite{teselucas, GouveiaQueiroz24, Pessoa24}).

In practice, applying Theorem \ref{TeoCyc1} involves computing the Jacobian matrix of the focus quantities with respect to the parameters under consideration and evaluating its rank at the point $s$ on the Bautin variety.

In this section, we will address the cyclicity of system \eqref{sist3Dchaotic} in a neighborhood of the singular point $E_1=(0,0,\frac{1}{d})$ in two scenarios. First, we will analyze the number of limit cycles that can arise simply by varying the parameters of system \eqref{sistE1norm} near the center condition. To do this, we will use the focus quantities obtained in the proof of Theorem \ref{Teo1} and a direct application of Theorem \ref{TeoCyc1}, using the linear parts of the focus quantities. Then, in a second step, we will perform a quadratic perturbation on system \eqref{sistE1cond6} with $d=1$, which has a center on the center manifold.
Naturally, if the perturbation parameters are zero, we have a center condition for the perturbed system. Hence, in this case, we will be dealing with quadratic perturbations of the original system, being able to generate more limit cycles using Theorem \ref{TeoCyc2}, by analyzing the higher orders of the focus quantities.

To be more precise, we first consider the following system
\begin{equation}
	\label{formnormPert1}
	\begin{array}{l}
		\dot{u} =\sigma u -v + P(u,v,w,\mu_1,\mu_2,\dots,\mu_n), \\
		\dot{v} = u + \sigma v+ Q(u,v,w,\mu_1,\mu_2,\dots,\mu_n), \\
		\dot{w} = \lambda w + R(u,v,w,\mu_1,\mu_2,\dots,\mu_n), \\
	\end{array}
\end{equation}
where $(\sigma,\lambda,\mu)\in\Lambda\subseteq\mathbb{R}\times\mathbb{R}^*\times\mathbb{R}^n$, with $\mu=(\mu_1,\mu_2,\dots\mu_n)$. Let $s$ be a point of the Bautin variety of system \eqref{formnormPert1}, with $\sigma=0$. In this case, computing the first $k$ focus quantities $L_1, L_2,\dots, L_k$, if the rank of the Jacobian matrix $J(L_1,L_2,\dots,L_k)=\left(\frac{\partial(L_1,L_2,\dots,L_k)}{\partial(\lambda,\mu_1,\mu_2,\dots\mu_n)}\right)_{\mu=s}$ is $k$, then we obtain $k-1$ limit cycles by a small variation of parameters $(\lambda,\mu)$, and together with a perturbation in the trace parameter $\sigma$, we obtain an additional limit cycle by the classical Hopf bifurcation. For this reason, in Theorem \ref{TeoCyc1}, we already include the total of $k$ limit cycles.

Next, to apply Theorem \ref{TeoCyc2}, we consider a system in the following form
\begin{equation}
	\label{formnormPert2}
	\begin{array}{l}
		\dot{u} = -v + P_c(u,v,w)+\displaystyle\sum_{j+k+l=2} a_{jkl}u^jv^kw^l, \\
		\dot{v} = u + Q_c(u,v,w)+\displaystyle\sum_{j+k+l=2} b_{jkl}u^jv^kw^l, \\
		\dot{w} = \lambda w + R_c(u,v,w)+\displaystyle\sum_{j+k+l=2} c_{jkl}u^jv^kw^l, \\
	\end{array}
\end{equation}
in which $\Lambda=(a_{jkl}, b_{jkl},c_{jkl}\;;\; j+k+l=2)\in\mathbb{R}^{18}$ are the perturbation parameters such that when $\Lambda=0$, we have a center at the origin on the center manifold. The analysis proceeds as follows: first, we calculate a sufficient number $N>k+l$ of focus quantities and compute the rank of the Jacobian matrix $J(L_1,\dots,L_N)=\left(\frac{\partial (L_1,\dots,L_N)}{\partial (a_{ijk},b_{ijk},c_{ijk})}\right)_{\Lambda=0}$. Let's say the rank is $k$. This means that we can obtain $k$ limit cycles by Theorem \ref{TeoCyc1} and the linear parts of the quantities $L_{k+1}, L_{k+2}, \dots, L_N$ are combinations of the linear parts of the first $k$ focus quantities, making it impossible to obtain more limit cycles from the linear parts of these quantities. Therefore, we must proceed to the analysis of the higher-order parts of these focus quantities. However, obtaining the polynomials $h_i$, according to Theorem \ref{TeoCyc2}, can be a computationally difficult task. Thus, to avoid performing this task directly and to obtain an equivalent result that is computationally easier to implement, we will make use of the following result, the proof of which can be found in \cite{teselucas,GouveiaTorre21,tesegouveia}.

%
% But, first, we emphasize here that the study of the number of limit cycles is clearly related to the study of the algebraic variety given by $S_L=\{ L_1=\cdots=L_N=0\}$, which is equivalent to the surjectivity of the mapping $L: \mathbb{R}^N \to \mathbb{R}^N$, $L(u)=(L_1(u),\dots, L_N(u))$,  near $u=0$, where $u=(u_1,\dots,u_N)$ (see \cite{teselucas,gouveiatesegouveia} for more details). The following result will be applied to simplify the focus quantities without modifying the algebraic variety $S_L$.

\begin{proposition}\label{CoroVariety}
	Let $f_i, g_i : \mathbb{R}^N \to \mathbb{R}$, $i = 1, \dots, k$, polynomial functions. Define the varieties $S_f = \{f_i = 0, i = 1, \dots, k\}$ and $S_g = \{g_i = 0, i = 1, \dots, k\}$. 
	If $g_i = f_i + \sum_{j<i} A_j^i f_j$ for $A_j^i : \mathbb{R}^N \to \mathbb{R}$ polynomial functions, then $S_g = S_f$.
\end{proposition}
First, let $k$ and $l$ be as in Theorem \ref{TeoCyc2}. Now, we make a linear change of variables to write 
\begin{equation*}
	L_i(u)= u_i+O(|u|^2), \qquad i=1,\dots,k,
\end{equation*}
and $u=(u_1,\dots,u_k,u_{k+1}, \dots,u_{k+l})$ the new parameters. Thus, the following $l$ focus quantities are written as
\begin{equation*}
	L_i(u)=\sum_{j=1}^{k}a_j^{i} u_j+O(|u|^2), \qquad i=k+1,\dots,k+l.
\end{equation*}
Defining new functions by 
\begin{equation*}
	\bar{L}_i(u)=L_i(u)-\sum_{j=1}^{k} a_l^{i} L_j(u), \qquad i=k+1,\dots,k+l,
\end{equation*}
we obtain, from Proposition \ref{CoroVariety}, that $S_L=\{L_1=\cdots =L_k=\bar{L}_{k+1}=\cdots=\bar{L}_{k+l}=0\}$, where $S_L=\{L_i=0,\;i=1,\dots,k+l\}$.
Furthermore, we have
\begin{equation*}
	\bar{L}_i(u)=h_i(u)+O(|u|^{m+1}),\qquad i=k+1,\dots,k+l,
\end{equation*}
where $h_i$ are homogeneous polynomials of degree $m\geq2$. At this point, in order to further simplify the variety $S_L$, we can take
\begin{equation*}
	\tilde{L}_i(u)=\bar{L}_i(u)-\Big(h_i(L_1,\dots,L_k,u_{k+1},\dots,u_{k+l})-h_i(0,\dots,0,u_{k+1},\dots,u_{k+l})\Big),\; i=k+1,\dots,k+l.
\end{equation*}
Hence, $S_L=\{L_1=\cdots= L_k=\tilde{L}_{k+1}=\cdots=\tilde{L}_{k+l}=0\}$ and the polynomials $\tilde{L}_i$, $i=k+1,\dots,k+l$, do not depend on the parameters $u_j$, $j=1,\dots,k$. From now on, we will work with the polynomials $L_1,\dots,L_k,\tilde{L}_{k+1},\dots,\tilde{L}_{k+l}$, but, to simplify the notation, we dropped the ``$\sim$" and return to just $L_i$, $i=1,\dots,k+l$.

Now, to verify the remaining hypotheses of Theorem \ref{TeoCyc2}, we look for a value $\eta$ of parameters which satisfies $h_i(\eta) = 0$, for $i = k+1, \dots, k+l-1$ and $h_{k+l}(\eta)\neq0$. To guarantee the transversality required in the theorem, the Jacobian matrix of $h_{k+1}, \dots, h_{k+l-1}$ (with respect to the variables $u_i$) must be non-singular at $\eta$. In the end, we obtain a total of $k+l$ limit cycles.

\subsection{Cyclicity without quadratic perturbation}
Let us consider the following system
\begin{equation}\label{sistE1normtrace}
	\begin{array}{l}
		\dot{u}=\sigma u+v+\dfrac{c d^2 (c d+1)}{c^2 d^2 k+k^3}u w+\dfrac{d \left(2 c^4 d^4+2 c^3 d^3+2 c^2 d^2 k^2+c^2 d^2+k^4\right)}{k^2 (c d+1) \left(c^2 d^2+k^2\right)}v w,\\ \\
		\dot{v}=-u+\sigma v-\dfrac{d (c d+1)}{c^2 d^2+k^2}u w-\dfrac{c d^2 (c d+1)}{k \left(c^2 d^2+k^2\right)}v w,\\ \\
		\dot{w}=-\dfrac{d^2}{k}w+\dfrac{d (c d+1)}{c^2 d^2+k^2}u v+\dfrac{c d^2 (c d+1)}{k \left(c^2 d^2+k^2\right)}v^2,
	\end{array}
\end{equation}
with parameters $\sigma,k,c$ and $d$. For $\sigma=0$,
in the proof of Theorem \ref{Teo1}, we see that this system has a center at the origin on the center manifold if and only if $k=1$ and $c=0$, leaving $d$ as a free parameter in the system.
\subsubsection{Proof of Theorem \ref{Teo4}}
Consider $L_1$, $L_2$, and $L_3$ the focus quantities of system \eqref{sistE1norm}, given in the proof of
Theorem \ref{Teo1}. Since the center condition of system \eqref{sistE1norm} is $k = 1$ and $c = 0$, the Jacobian matrix of these quantities with respect to the parameters $(k, c, d)$, evaluated at the point $s = (1, 0, d)$ of the associated Bautin variety, is
\begin{equation*}
	\left(\dfrac{\partial(L_1,L_2,L_3)}{\partial(k,c,d)}\right)_{(k,c,d)=s}=\begin{pmatrix}
		6 & 16d+46 & 204d+36 \\
		80 & P_{22}(d) & P_{23}(d) \\
		d^5 & P_{32}(d) & P_{33}(d)
	\end{pmatrix},
\end{equation*}
where 
\begin{flalign*}
	P_{22}(d) &= 256 d^{8}+248 d^{7}-32 d^{5}+15448 d^{4}+12288 d^{3}+428 d+4048, & \\
	P_{23}(d) &= 378 d^{6}+21816 d^{4}-7776 d^{3}+139176 d^{2}+30228 d+289584, &
\end{flalign*}

\begin{flalign*}
	P_{32}(d) &= 1697280 d^{31}+2350080 d^{30}+1934464 d^{29}+844288 d^{28}+32059456 d^{27}+ & \\
	& \quad 43939072 d^{26}+36594208 d^{25}+16221440 d^{24}-31883648 d^{23}-54237952 d^{22} & \\
	& \quad -43917536 d^{21}-22140160 d^{20}-1765557696 d^{19}-2449146880 d^{18}- & \\
	& \quad 2025351296 d^{17}-923707392 d^{16}-6577764864 d^{15}-8741859328 d^{14}- & \\
	& \quad 7237094912 d^{13}-3307204608 d^{12}-4621947904 d^{11}-4961091584 d^{10}- & \\
	& \quad 3855288320 d^{9}-1644658688 d^{8}-75497472 d^{6}-509648896 d^{5}, &
\end{flalign*}

\begin{flalign*}
	P_{33}(d) &= 470239992 d^{30}+463993920 d^{29}+276975288 d^{28}+89558784 d^{27}+ & \\
	& \quad 7665271704 d^{26}+7495690176 d^{25}+4493995800 d^{24}+1466284032 d^{23}- & \\
	& \quad 7968309888 d^{22}-8377737984 d^{21}-5105395008 d^{20}-1894682880 d^{19}- & \\
	& \quad 301074259008 d^{18}-291380405760 d^{17}-172641146112 d^{16}- & \\
	& \quad 57169649664 d^{15}-855101836800 d^{14}-802525151232 d^{13}- & \\
	& \quad 468463816704 d^{12}-152505815040 d^{11}-372904080384 d^{10}- & \\
	& \quad 312222228480 d^{9}-94202910720 d^{8}-18795110400 d^{7}-10653081600 d^{6}- & \\
	& \quad 76508946432 d^{5}-31431352320 d^{4}. &
\end{flalign*}
Since it has rank $3$, by applying Theorem \ref{TeoCyc1} to system \eqref{sistE1normtrace}, we conclude that there exist variations in the parameters of system \eqref{sist3Dchaotic} for which this system has at least three limit cycles. \qed

\subsection{Cyclicity with quadratic perturbation}
Consider system \eqref{sistE1cond6} with the center conditions $k=1$ and $c=0$. We also assume $d=1$. Thus, performing a quadratic perturbation as in \eqref{formnormPert2}, we have the following
 \begin{equation}\label{sistE1cond6Pert}
 	\begin{array}{l}
 		\dot{u}=v+ v w+\displaystyle\sum_{j+k+l=2} a_{jkl}u^jv^kw^l,\\
 		\dot{v}=-u- u w+\displaystyle\sum_{j+k+l=2} b_{jkl}u^jv^kw^l, \\
 		\dot{w}=-w+ u v+\displaystyle\sum_{j+k+l=2} c_{jkl}u^jv^kw^l,
 	\end{array}
 \end{equation}
with $\Lambda=(a_{jkl}, b_{jkl},c_{jkl}\;;\; j+k+l=2)\in\mathbb{R}^{18}$ the perturbation parameters. Note that, for $\Lambda=0$, we have a center at the origin on the center manifold. 

\subsubsection{Proof of Theorem \ref{Teo5}} Firstly, denoting by $L_k^j$ the homogeneous part of order $j$ of the $k$--th focus quantities $L_k$, we compute the linear part of the first 9 focus quantities for system \eqref{sistE1cond6Pert}, that is,
\begin{align*}
	L_1^1 &= \frac{1}{20} (a_{011} + 2\,a_{101} - 2\,b_{011} + b_{101}); \\
	L_2^1 &= \frac{1}{40} (-a_{101} - b_{011}); \\
	L_3^1 &= \frac{281\,a_{011} + 342\,a_{101} - 342\,b_{011} + 281\,b_{101}}{136000}; \\
	L_4^1 &= -\frac{281 (a_{101} + b_{011})}{272000}; \\
	L_5^1 &= \frac{2324157\,a_{011} + 2420774\,a_{101} - 2420774\,b_{011} + 2324157\,b_{101}}{17108800000}; \\
	L_6^1 &= -\frac{2324157 (a_{101} + b_{011})}{34217600000}; \\
	L_7^1 &= \frac{18296103569\,a_{011} + 17579350678\,a_{101} - 17579350678\,b_{011} + 18296103569\,b_{101}}{1721829632000000}; \\
	L_8^1 &= -\frac{18296103569 (a_{101} + b_{011})}{3443659264000000}; \\
	L_9^1 &= \frac{1}{1422019490987264000000000} (1295884288642940083\,a_{011} + 1183999528745548106\,a_{101} \\
	& \quad - 1183999528745548106\,b_{011} + 1295884288642940083\,b_{101}).
\end{align*}
We have that \begin{equation*}
	\textrm{rank}J(L_1,L_2,\dots,L_9)_{\Lambda=0}=\textrm{rank}J(L_1,L_2,L_3)_{\Lambda=0}.
\end{equation*}
Hence, by Theorem \ref{TeoCyc1}, is possible to obtain 3 limit cycles, by perturbation on $\Lambda$, from the center at the origin. Now, after a suitable change of variables in the parameters $(a_{011},a_{101},b_{011})$, and some simplifications, as described at the beginning of this section, we can write 
\begin{equation*}
	L_i=u_i+O(\vert u \vert^2),\quad i=1,2,3, \qquad \textrm{and}\qquad
	L_j=h_j+O(\vert u \vert^3),\quad j=4,5,
\end{equation*}
where $u=(a_{200}, a_{110}, a_{020}, a_{002}, b_{200}, b_{110}, b_{101}, b_{020}, 
b_{002}, c_{200}, c_{110}, c_{101}, c_{020}, c_{011}, c_{002})$ and 
\begin{align*}
	h_4 &= \frac{147691\,a_{002}^2}{106080000} + \frac{24061\,a_{002} a_{020}}{35360000} + \frac{19453\,a_{020}^2}{3264000} - \frac{229649\,a_{002} a_{110}}{106080000} - \frac{177\,a_{020} a_{110}}{544000} \\
	&\quad + \frac{2699\,a_{110}^2}{3264000} - \frac{274571\,a_{002} a_{200}}{35360000} + \frac{433\,a_{020} a_{200}}{32000} + \frac{3371\,a_{110} a_{200}}{544000} + \frac{22047\,a_{200}^2}{1088000} \\
	&\quad - \frac{27639\,a_{002} b_{002}}{17680000} - \frac{106841\,a_{020} b_{002}}{21216000} + \frac{12607\,a_{110} b_{002}}{7072000} - \frac{36311\,a_{200} b_{002}}{7072000} - \frac{147691\,b_{002}^2}{106080000} \\
	&\quad + \frac{36311\,a_{002} b_{020}}{7072000} - \frac{4571\,a_{020} b_{020}}{544000} - \frac{19\,a_{110} b_{020}}{4352} - \frac{981\,a_{200} b_{020}}{32000} - \frac{274571\,b_{002} b_{020}}{35360000} - \frac{22047\,b_{020}^2}{1088000} \\
	&\quad + \frac{12607\,a_{002} b_{110}}{7072000} + \frac{8377\,a_{020} b_{110}}{1632000} - \frac{1053\,a_{110} b_{110}}{544000} + \frac{19\,a_{200} b_{110}}{4352} + \frac{229649\,b_{002} b_{110}}{106080000} \\
	&\quad + \frac{3371\,b_{020} b_{110}}{544000} - \frac{2699\,b_{110}^2}{3264000} + \frac{106841\,a_{002} b_{200}}{21216000} + \frac{699\,a_{020} b_{200}}{544000} - \frac{8377\,a_{110} b_{200}}{1632000} - \frac{4571\,a_{200} b_{200}}{544000} \\
	&\quad + \frac{24061\,b_{002} b_{200}}{35360000} - \frac{433\,b_{020} b_{200}}{32000} - \frac{177\,b_{110} b_{200}}{544000} - \frac{19453\,b_{200}^2}{3264000} - \frac{27241\,a_{002} c_{011}}{35360000} - \frac{22297\,a_{020} c_{011}}{14144000} \\
	&\quad + \frac{8089\,a_{110} c_{011}}{14144000} + \frac{3759\,a_{200} c_{011}}{14144000} + \frac{13287\,b_{002} c_{011}}{8840000} + \frac{60577\,b_{020} c_{011}}{14144000} - \frac{22297\,b_{110} c_{011}}{14144000}\\
	&\quad  + \frac{8089\,b_{200} c_{011}}{14144000} + \frac{13287\,a_{002} c_{101}}{8840000} - \frac{8089\,a_{020} c_{101}}{14144000} - \frac{22297\,a_{110} c_{101}}{14144000} - \frac{60577\,a_{200} c_{101}}{14144000}  \\
	&\quad + \frac{27241\,b_{002} c_{101}}{35360000} + \frac{3759\,b_{020} c_{101}}{14144000} - \frac{8089\,b_{110} c_{101}}{14144000} - \frac{22297\,b_{200} c_{101}}{14144000}; 
	\end{align*}
	 
	 \begin{align*}
	 	h_5 &= -\frac{1144493\,a_{002}^2}{5233280000} + \frac{62982333\,a_{002} a_{020}}{71957600000} + \frac{107984027\,a_{020}^2}{172698240000} + \frac{9209861\,a_{002} a_{110}}{39249600000} - \frac{56388799\,a_{020} a_{110}}{57566080000}\\
	 	&\quad  - \frac{1000277\,a_{110}^2}{15699840000} + \frac{287805833\,a_{002} a_{200}}{107936400000} - \frac{379313\,a_{020} a_{200}}{1569984000} - \frac{3319067\,a_{110} a_{200}}{1381585920}- \frac{220968679\,a_{200}^2}{57566080000}\\
	 	&\quad   - \frac{337319707\,a_{020} b_{002}}{863491200000} + \frac{1180287\,a_{110} b_{002}}{287830400000}- \frac{193767471\,a_{200} b_{002}}{287830400000} - \frac{1144493\,b_{002}^2}{5233280000} \\
	 	&\quad  - \frac{193767471\,a_{002} b_{020}}{287830400000} + \frac{6270373\,a_{020} b_{020}}{43174560000} + \frac{43572101\,a_{110} b_{020}}{43174560000} - \frac{287805833\,b_{002} b_{020}}{107936400000} \\
	 	&\quad - \frac{220968679\,b_{020}^2}{57566080000} - \frac{1180287\,a_{002} b_{110}}{287830400000} + \frac{4849049\,a_{020} b_{110}}{8634912000} + \frac{43572101\,a_{200} b_{110}}{43174560000}+ \frac{9209861\,b_{002} b_{110}}{39249600000}\\
	 	&\quad   + \frac{3319067\,b_{020} b_{110}}{1381585920} - \frac{1000277\,b_{110}^2}{15699840000} - \frac{337319707\,a_{002} b_{200}}{863491200000} + \frac{4849049\,a_{110} b_{200}}{8634912000} - \frac{6270373\,a_{200} b_{200}}{43174560000}\\
	 	&\quad - \frac{62982333\,b_{002} b_{200}}{71957600000} - \frac{379313\,b_{020} b_{200}}{1569984000} + \frac{56388799\,b_{110} b_{200}}{57566080000} + \frac{107984027\,b_{200}^2}{172698240000} - \frac{93582813\,a_{002} c_{011}}{575660800000}\\
	 	&\quad  - \frac{38776721\,a_{020} c_{011}}{115132160000} + \frac{260549\,a_{110} c_{011}}{2214080000} + \frac{10818981\,a_{200} c_{011}}{115132160000} + \frac{139511491\,b_{002} c_{011}}{575660800000} \\
	 	&\quad  + \frac{18348333\,b_{020} c_{011}}{28783040000}- \frac{38776721\,b_{110} c_{011}}{115132160000} + \frac{260549\,b_{200} c_{011}}{2214080000}- \frac{139511491\,a_{002} c_{101}}{575660800000}  \\
	 	&\quad + \frac{260549\,a_{020} c_{101}}{2214080000} + \frac{38776721\,a_{110} c_{101}}{115132160000} + \frac{18348333\,a_{200} c_{101}}{28783040000}- \frac{93582813\,b_{002} c_{101}}{575660800000} \\
	 	&\quad  - \frac{10818981\,b_{020} c_{101}}{115132160000} + \frac{260549\,b_{110} c_{101}}{2214080000} + \frac{38776721\,b_{200} c_{101}}{115132160000}.
	 \end{align*}

Taking the line $\eta$, from the statement of Theorem \ref{TeoCyc2}, as such $\{a_{002} = a_{110} = a_{020} = a_{200} = b_{002} = b_{020} = b_{110} = c_{011} = 0, c_{101} = - \frac{252889 b_{200}}{66891}\}$, we have $L_4(\eta)=0$ and $L_5(\eta)=-\frac{4990766496931\, b_{200}^2}{7701305314560000}\neq 0$. This implies, by Theorem \ref{TeoCyc2}, that it is possible to obtain 2 additional limit cycles. In other words, through quadratic perturbations, at least 5 limit cycles can bifurcate from the origin of system \eqref{sistE1cond6}. \qed 

\appendix
\section{}\label{apendice}
\begin{align*}
	L_3&=d^5 \bigg(816 c^9 d^{26}-408 c^7 d^{26}+2448 c^8 d^{25}+204 k^2 c^6 d^{25}-1020 c^6 d^{25}+776 k^2 c^7 d^{24}+3672 c^7 d^{24}\\
	&\quad+496 k^2 c^5 d^{24}-1224 c^5 d^{24}+1104 k^2 c^6 d^{23}+3264 c^6 d^{23}-44 k^4 c^4 d^{23}+1036 k^2 c^4 d^{23}\\
	&\quad-816 c^4 d^{23}+15152 k^2 c^9 d^{22}-6872 k^2 c^7 d^{22}+452 k^4 c^5 d^{22}+544 k^2 c^5 d^{22}+1836 c^5 d^{22}\\
	&\quad-362 k^4 c^3 d^{22}+780 k^2 c^3 d^{22}-306 c^3 d^{22}+45456 k^2 c^8 d^{21}+3436 k^4 c^6 d^{21}-17180 k^2 c^6 d^{21}\\
	&\quad+692 k^4 c^4 d^{21}-344 k^2 c^4 d^{21}+612 c^4 d^{21}+35 k^6 c^2 d^{21}-105 k^4 c^2 d^{21}+185 k^2 c^2 d^{21}-51 c^2 d^{21}\\
	&\quad+14856 k^4 c^7 d^{20}+68184 k^2 c^7 d^{20}+12976 k^4 c^5 d^{20}-20616 k^2 c^5 d^{20}+50 k^6 c^3 d^{20}+642 k^4 c^3 d^{20}\\
	&\quad-506 k^2 c^3 d^{20}+102 c^3 d^{20}-38 k^6 c d^{20}+76 k^4 c d^{20}-6 k^2 c d^{20}+k^8 d^{19}-19 k^6 d^{19}\\
	&\quad+21840 k^4 c^6 d^{19}+60608 k^2 c^6 d^{19}+19 k^4 d^{19}-3052 k^6 c^4 d^{19}+29004 k^4 c^4 d^{19}-13744 k^2 c^4 d^{19}\\
	&\quad-k^2 d^{19}+148 k^6 c^2 d^{19}+296 k^4 c^2 d^{19}-252 k^2 c^2 d^{19}-21760 k^4 c^9 d^{18}+32064 k^4 c^7 d^{18}\\
	&\quad+6980 k^6 c^5 d^{18}+6304 k^4 c^5 d^{18}+34092 k^2 c^5 d^{18}-9818 k^6 c^3 d^{18}+23020 k^4 c^3 d^{18}-5154 k^2 c^3 d^{18}\\
	&\quad-42 k^8 c d^{18}+74 k^6 c d^{18}+74 k^4 c d^{18}-42 k^2 c d^{18}-65280 k^4 c^8 d^{17}-16032 k^6 c^6 d^{17}\\
	&\quad+80160 k^4 c^6 d^{17}+13748 k^6 c^4 d^{17}-16216 k^4 c^4 d^{17}+11364 k^2 c^4 d^{17}+139 k^8 c^2 d^{17}\\
	&\quad-2785 k^6 c^2 d^{17}+6385 k^4 c^2 d^{17}-859 k^2 c^2 d^{17}-34944 k^6 c^7 d^{16}-97920 k^4 c^7 d^{16}+99776 k^6 c^5 d^{16}\\
	&\quad+96192 k^4 c^5 d^{16}-1230 k^8 c^3 d^{16}+11778 k^6 c^3 d^{16}-18554 k^4 c^3 d^{16}+1894 k^2 c^3 d^{16}-1062 k^8 c d^{16}\\
	&\quad+2124 k^6 c d^{16}+378 k^4 c d^{16}-63 k^{10} d^{15}-531 k^8 d^{15}+531 k^6 d^{15}-72192 k^6 c^6 d^{15}-87040 k^4 c^6 d^{15}\\
	&\quad+63 k^4 d^{15}-65920 k^8 c^4 d^{15}+265472 k^6 c^4 d^{15}+64128 k^4 c^4 d^{15}+2452 k^8 c^2 d^{15}+4904 k^6 c^2 d^{15}\\
	&\quad-8700 k^4 c^2 d^{15}-859904 k^6 c^9 d^{14}+515456 k^6 c^7 d^{14}-31296 k^8 c^5 d^{14}-133760 k^6 c^5 d^{14}\\
	&\quad-48960 k^4 c^5 d^{14}-93584 k^8 c^3 d^{14}+224416 k^6 c^3 d^{14}+24048 k^4 c^3 d^{14}-1450 k^{10} c d^{14}+1226 k^8 c d^{14}\\
	&\quad+1226 k^6 c d^{14}-1450 k^4 c d^{14}-2579712 k^6 c^8 d^{13}-257728 k^8 c^6 d^{13}+1288640 k^6 c^6 d^{13}\\
	&\quad+320 k^8 c^4 d^{13}-158080 k^6 c^4 d^{13}-16320 k^4 c^4 d^{13}-11112 k^{10} c^2 d^{13}-31624 k^8 c^2 d^{13}\\
	&\quad+67144 k^6 c^2 d^{13}+4008 k^4 c^2 d^{13}-977536 k^8 c^7 d^{12}-3869568 k^6 c^7 d^{12}+301312 k^8 c^5 d^{12}\\
	&\quad+1546368 k^6 c^5 d^{12}-24736 k^{10} c^3 d^{12}+2720 k^8 c^3 d^{12}-115040 k^6 c^3 d^{12}-2720 k^4 c^3 d^{12}\\
	&\quad-7584 k^{10} c d^{12}+15168 k^8 c d^{12}+6624 k^6 c d^{12}-1104 k^{12} d^{11}-3792 k^{10} d^{11}+3792 k^8 d^{11}\\
	&\quad+1104 k^6 d^{11}-1642752 k^8 c^6 d^{11}-3439616 k^6 c^6 d^{11}-408384 k^{10} c^4 d^{11}+1011008 k^8 c^4 d^{11}\\
	&\quad+1030912 k^6 c^4 d^{11}+1280 k^{10} c^2 d^{11}+2560 k^8 c^2 d^{11}-45312 k^6 c^2 d^{11}-3018752 k^8 c^9 d^{10}\\
	&\quad+1541120 k^8 c^7 d^{10}-543040 k^{10} c^5 d^{10}-1568256 k^8 c^5 d^{10}-1934784 k^6 c^5 d^{10}-386656 k^{10} c^3 d^{10}\\
	&\quad+893248 k^8 c^3 d^{10}+386592 k^6 c^3 d^{10}-7552 k^{12} c d^{10}+640 k^{10} c d^{10}+640 k^8 c d^{10}-7552 k^6 c d^{10}\\
	&\quad-9056256 k^8 c^8 d^9-770560 k^{10} c^6 d^9+3852800 k^8 c^6 d^9-413760 k^{10} c^4 d^9-828544 k^8 c^4 d^9\\
	&\quad-644928 k^6 c^4 d^9-86192 k^{12} c^2 d^9-155632 k^{10} c^2 d^9+264432 k^8 c^2 d^9+64432 k^6 c^2 d^9 \\
	&\quad-3434496 k^{10} c^7 d^8-13584384 k^8 c^7 d^8+345088 k^{10} c^5 d^8+4623360 k^8 c^5 d^8-155040 k^{12} c^3 d^8\\
	&\quad-274592 k^{10} c^3 d^8-309216 k^8 c^3 d^8-107488 k^6 c^3 d^8-18848 k^{12} c d^8+37696 k^{10} c d^8+24672 k^8 c d^8\\
	&\quad-4112 k^{14} d^7-9424 k^{12} d^7+9424 k^{10} d^7+4112 k^8 d^7-5775360 k^{10} c^6 d^7-12075008 k^8 c^6 d^7\\
	&\quad-943104 k^{12} c^4 d^7+1633280 k^{10} c^4 d^7+3082240 k^8 c^4 d^{7}-33856 k^{12} c^2 d^7-67712 k^{10} c^2 d^7\\
	&\quad-61248 k^8 c^2 d^7-1572864 k^{10} c^9 d^6+131072 k^{10} c^7 d^6-1758208 k^{12} c^5 d^6-4950016 k^{10} c^5 d^6\\
	&\quad-6792192 k^8 c^5 d^6-699648 k^{12} c^3 d^6+1485312 k^{10} c^3 d^6+1155840 k^8 c^3 d^6-10208 k^{14} c d^6\\
	&\quad-16928 k^{12} c d^6-16928 k^{10} c d^6-10208 k^8 c d^6-4718592 k^{10} c^8 d^5-65536 k^{12} c^6 d^5+327680 k^{10} c^6 d^5\\
	&\quad-1571840 k^{12} c^4 d^5-1783808 k^{10} c^4 d^5-2264064 k^8 c^4 d^5-208000 k^{14} c^2 d^5-325248 k^{12} c^2 d^5\\
	&\quad+402048 k^{10} c^2 d^5+192640 k^8 c^2 d^5-1572864 k^{12} c^7 d^4-7077888 k^{10} c^7 d^4+131072 k^{12} c^5 d^4\\
	&\quad+393216 k^{10} c^5 d^4-328192 k^{14} c^3 d^4-1007104 k^{12} c^3 d^4-97792 k^{10} c^3 d^4-377344 k^8 c^3 d^4\\
	&\quad-12288 k^{14} c d^4+24576 k^{12} c d^4+18432 k^{10} c d^4-3072 k^{16} d^3-6144 k^{14} d^3+6144 k^{12} d^3+3072 k^{10} d^3\\
	&\quad-2359296 k^{12} c^6 d^3-6291456 k^{10} c^6 d^3-131072 k^{14} c^4 d^3+393216 k^{12} c^4 d^3+262144 k^{10} c^4 d^3\\
	&\quad-110592 k^{14} c^2 d^3-221184 k^{12} c^2 d^3+208896 k^{10} c^2 d^3-786432 k^{14} c^5 d^2-1966080 k^{12} c^5 d^2\\
	&\quad-3538944 k^{10} c^5 d^2-163840 k^{14} c^3 d^2+327680 k^{12} c^3 d^2+98304 k^{10} c^3 d^2+34816 k^{16} c d^2\\
	&\quad-55296 k^{14} c d^2-55296 k^{12} c d^2+34816 k^{10} c d^2-393216 k^{14} c^4 d-786432 k^{12} c^4 d-1179648 k^{10} c^4 d\\
	&\quad-16384 k^{16} c^2 d-81920 k^{14} c^2 d+81920 k^{12} c^2 d+16384 k^{10} c^2 d-196608 k^{16} c^3-196608 k^{14} c^3\\
	&\quad-196608 k^{12} c^3-196608 k^{10} c^3\bigg).
\end{align*}

\section*{Author contributions}

C.P. proposed the research problem and the methodology to address it, contributed to the conceptualization, data curation, and resources, and provided input during the analysis and final revision of the manuscript. V.G. conducted the literature review, prepared the figures, and wrote the the first draft. Both authors contributed jointly to the mathematical development and proofs of the results, critically reviewed the manuscript’s content, and approved the final version.

\section*{Funding }

The first author is supported by CAPES grant 88887.952598/2024-00. The second author is partially supported by S\~ao Paulo Research Foundation (FAPESP) grants 24/16727-1 and 24/15612-6.

 \section*{Data availability}
The datasets generated during and/or analyzed during the current study are available from the corresponding author on reasonable request.

\section*{Declarations}
 
 \section*{Conflict of interest}
 The authors have no conflict of interest.
 
  \section*{Ethical standard} 
  All procedures performed in studies
involving human participants were in accordance with the ethical standards of the institutional and/or national research
committee and with the 1964 Helsinki Declaration and its later
amendments or comparable ethical standards.

 \section*{Human and animal rights} 
 This article does not contain any
studies with animals performed by any of the authors.

 \section*{Informed consent} Informed consent was obtained from all
individual participants included in the study.

%\bibliographystyle{siam}
%\bibliography{Ref_3.bib}

\end{document}